\documentclass[12 pt]{article}
\usepackage[pagebackref]{hyperref}
\hypersetup{
colorlinks = true,
linkcolor = blue,
citecolor = blue}
\usepackage[utf8]{inputenc}
\usepackage[T1]{fontenc}

\usepackage[]{amsfonts}
\usepackage{amssymb}
\usepackage{amsmath,amsthm}
\def\couleur(#1 #2 #3)
	{
\special{ps::[end]}}

\def\bx#1{\setbox1=\hbox{\kern3pt{#1}\kern3pt}			% Make a box. Close it by "}"
 \dimen1=\ht1 \advance\dimen1 by 3pt \dimen2=\dp1 \advance\dimen2 by 3pt
 \setbox1=\hbox{\vrule height\dimen1 depth\dimen2\box1\vrule}%
 \setbox1=\vbox{\hrule\box1\hrule}%
 \advance\dimen1 by .4pt \ht1=\dimen1
 \advance\dimen2 by .4pt \dp1=\dimen2 \box1\relax}

\def\wbb#1{\kern#1em}
\def\vci{\vrule  width.02em height1.47ex depth-.0ex}		% le 1 en blackboard
\def\11{{\rm\wbb{.2}\vci\wbb{-.37}1}}

\def\underset#1#2{\mathrel{\mathop{\kern0pt #2}\limits_{#1}}}

\def\overset#1#2{\mathrel{\mathop{\kern0pt #2}\limits^{#1}}}

\newtheorem{thm}{Theorem}[section]
\newtheorem{lem}[thm]{Lemma}
\newtheorem{prop}[thm]{Proposition}
\newtheorem{defin}[thm]{Definition}
\newtheorem{rem}[thm]{Remark}

\begin{document}

\title{On the linear extension property for interpolating sequences.}

\author{Eric Amar}

\date{}
\maketitle
 \renewcommand{\abstractname}{Abstract}

\begin{abstract}
Let $S$ be a sequence of points in $\Omega ,$ where  $\Omega
 $ is the unit ball or the unit polydisc in ${\mathbb{C}}^{n}.$
 Denote $H^{p}$($\Omega $) the Hardy space of $\Omega .$ Suppose
 that $S$ is $H^{p}$ interpolating with $p\geq 2.$ Then $S$ has
 the bounded linear extension property. The same is true for
 the Bergman spaces of the ball by use of the "Subordination
 Lemma". The point of view used here is the vectorial one: Hilbertian
 and Besselian basis.\par 

\end{abstract}
\ \par 
\ \par 
\ \par 
\ \par 

\tableofcontents
\ \par 

\section{Introduction.}
\quad Let ${\mathbb{D}}$ be the unit disc in the complex plane and
 $S$ be a sequence of points in ${\mathbb{D}}.$ Denote $H^{\infty
 }$ the set of \textbf{bounded holomorphic functions} on ${\mathbb{D}}.$
 We shall say that $S$ is a $H^{\infty }$ \textbf{interpolating
 sequence} if\ \par 
\quad \quad \quad $\forall \lambda \in ???^{\infty }(S),\exists f\in H^{\infty
 }::\forall a\in S,\ f(a)=\lambda _{a}.$\ \par 
These interpolating sequences where characterized by L. Carleson~\cite{CarlInt58}
 and by H. Shapiro and A. Shieds~\cite{ShapShields61} for the
 Hardy spaces $H^{p}({\mathbb{D}})$ by the same condition:\ \par 
\quad \quad \quad $\displaystyle \forall b\in S,\ \prod_{a\in S,\ a\neq b}{\left\vert{\frac{a-b}{1-\bar
 ab}}\right\vert }\geq \delta >0.$\ \par 
\quad In several variables for the unit ball $\Omega ={\mathbb{B}}\subset
 {\mathbb{C}}^{n}$ or for the unit polydisc $\Omega ={\mathbb{D}}^{n}\subset
 {\mathbb{C}}^{n},$ this characterisation is still an open question,
 even for the Hilbert space $H^{2}(\Omega ).$\ \par 
\quad The aim is to study interpolating sequences in the Hardy spaces
 in the ball and in the polydisc in ${\mathbb{C}}^{n}.$\ \par 
\quad We shall follow here the \emph{vectorial} point of view, as in
 ~\cite{ChalChev18}, ~\cite{AmarThesis77}, ~\cite{HardySob18}
 or in ~\cite{AglMCar02}, instead of the \emph{functional} one,
 as done classically in~\cite{ShapShields61}.\ \par 
\quad For $p\geq 1$ we note $p'$ its conjugate $\frac{1}{p'}+\frac{1}{p}=1.$
 For a Banach space ${\mathcal{B}},$ call ${\mathcal{B}}'$ its
 dual space.\ \par 

\begin{defin}
Let $p>1$ and ${\mathcal{S}}:=\lbrace e_{a}\rbrace _{a\in S}$
 a sequence of vectors of norm one in the Banach space ${\mathcal{B}}'.$
 Then ${\mathcal{S}}$ is $p$-{\bf Carleson} if it is $p'$-{\bf
 hilbertian} i.e.:\par 
\quad \quad \quad $(H_{p'})\ \ \ \ \ \ \exists C>0,\ \forall \lambda \in ???^{p'}(S),\
 {\left\Vert{\sum_{a\in S}{\lambda _{a}e_{a}}}\right\Vert}_{{\mathcal{B}}'}\leq
 C{\left\Vert{\lambda }\right\Vert}_{???^{p'}(S)}.$\par 
${\mathcal{S}}$ is $p$-{\bf interpolating} if it is $p'$-{\bf
 besselian} i.e.:\par 
\quad \quad \quad $(B_{p'})\ \ \ \ \ \ \exists C>0,\ \forall \lambda \in ???^{p'}(S),\
 {\left\Vert{\sum_{a\in S}{\lambda _{a}e_{a}}}\right\Vert}_{{\mathcal{B}}'}\geq
 c{\left\Vert{\lambda }\right\Vert}_{???^{p'}(S)}.$\par 
It is a $p'$-{\bf Riesz} sequence if ${\mathcal{S}}$ is $p'$-hilbertian
 and $p'$-besselian.
\end{defin}
\quad Recall a result of Babenko~\cite{Babenko47}: there exists $\lbrace
 e_{a}\rbrace _{a\in S}$ a basis of unit vectors in the Hilbert
 space ${\mathcal{H}}$ which is $2$-\emph{hilbertian} but $\lbrace
 e_{a}\rbrace _{a\in S}$ is not a Riesz basis.\ \par 
\quad The same way there exists $\lbrace e_{a}\rbrace _{a\in S}$ a
 basis of unit vectors in the Hilbert space ${\mathcal{H}}$ which
 is $2$-\emph{besselian} but $\lbrace e_{a}\rbrace _{a\in S}$
 is not a Riesz basis.\ \par 
\quad To have a elementary proof of Babenko results, together with
 an extension of Babenko examples to $???^{p}$ spaces, see the
 work by I. Chalendar, B. Chevreau and E. A.~\cite{ChalChev18}.\ \par 
\quad We shall get as a by-product of the results here that we have
 Hardy spaces of the ball or of the polydisc in ${\mathbb{C}}^{2}$
 which are $p$-hilbertian but not $p$-besselian, giving new Babenko
 examples.\ \par 
\ \par 
\quad We shall apply this definition to the case of ${\mathcal{B}}=H^{p}(\Omega
 ),$ the Hardy space of $\Omega ,$ with $\Omega $ the ball or
 the polydisc in ${\mathbb{C}}^{n},$ and with $e_{a}=k_{a,p'}\in
 H^{p'}=(H^{p})^{*},$ the normalised reproducing kernel of the
 point $\displaystyle a\in \Omega .$ The precise definitions
 are in the next section.\ \par 
\quad Define the restriction operator $R_{p}:H^{p}(\Omega )\rightarrow
 ???^{0}(S),$ with $???^{0}(S)$ the set of complex valued sequences, by:\ \par 
\quad \quad \quad $\forall f\in H^{p}(\Omega ),\ R_{p}f:=\lbrace {\left\langle{f,k_{a,p'}}\right\rangle}\rbrace
 _{a\in S}\in ???^{0}(S).$\ \par 
\quad As we shall see in Lemma~\ref{pH1}, we have that $S$ is $H^{p}(\Omega
 )$ interpolating iff we have that $R_{p}(H^{p}(\Omega ))\supset
 ???^{p}(S),$ hence, for any sequence $\lambda \in ???^{p}(S),$
 there is a function $f\in H^{p}$ such that $\forall a\in S,\
 {\left\langle{f,k_{a,p'}}\right\rangle}=\lambda _{a}.$ So the
 vectorial definition of $p$-interpolating sequences is equivalent
 to the functional usual one.\ \par 
\quad Now we can set\ \par 

\begin{defin}
Let $p\geq 1$ and $S\subset \Omega $ be a sequence of points
 in $\Omega .$ We say that the sequence $S$ has the {\bf linear
 extension property, LEP,} if there is a bounded linear operator
 $E:\ell ^{p}(S)\rightarrow H^{p}$ such that for any $\lambda
 \in \ell ^{p}(S),\ E(\lambda )$ interpolates the sequence $\lambda
 $ in $H^{p}(\Omega ).$  I.e. $\displaystyle \forall a\in S,\
 {\left\langle{E(\lambda ),k_{a,p'}}\right\rangle}=\lambda _{a}.$\par 
\quad Its range $E_{S}^{p}$ is the subspace $E(???^{p}(S))$ of $H^{p}.$
\end{defin}
\quad P. Beurling~\cite{PBeurling62} proved that the $H^{\infty }({\mathbb{D}})$
 interpolating sequences in the unit disc have always the LEP.\ \par 
Using a very nice method due to S. Drury~\cite{Drury70}, A. Bernard~\cite{Bernard71}
 proved the same for interpolating sequences in uniform algebras,
 hence for $H^{\infty }({\mathbb{D}}^{n})$ and $H^{\infty }({\mathbb{B}}).$\
 \par 
\quad For $H^{p}({\mathbb{D}})$ interpolating sequences, $1\leq p<\infty
 ,$ I proved in~\cite{amExt83} that they have the LEP, by use
 of $\bar \partial $ methods. In~\cite{SchuSeip98} this result
 is reproved by a different method.\ \par 
\quad This question was open for a while in the several variables case. \ \par 
\ \par 
\quad Using the vectorial point of view, we shall prove:\ \par 

\begin{thm}
Let $S$ be a sequence of points in $\Omega .$ Suppose that $S$
 is $H^{p}(\Omega )$ interpolating with $p\geq 2.$ Then $S$ has the LEP.
\end{thm}
\quad We also study \textbf{strictly }$H^{p}$ interpolating sequences
 and the case of interpolating sequences for \textbf{weighted
 Bergman spaces.}\ \par 

\section{Definitions, notation.}
\quad Recall the definition of Hardy spaces in the polydisc.\ \par 

\begin{defin}
The {\bf Hardy space} $H^{p}({\mathbb{D}}^{n})$ is the set of
 holomorphic functions $f$ in ${\mathbb{D}}^{n}$ such that, with
 $e^{i\theta }:=e^{i\theta _{1}}{\times}\cdot \cdot \cdot {\times}e^{i\theta
 _{n}}$ and $d\theta :=d\theta _{1}\cdot \cdot \cdot d\theta
 _{n}$ the Lebesgue measure on ${\mathbb{T}}^{n}$:\par 
\quad \quad \quad $\displaystyle {\left\Vert{f}\right\Vert}_{H^{p}}^{p}:=\sup _{r<1}\int_{{\mathbb{T}}^{n}}{\left\vert{f(re^{i\theta
 })}\right\vert ^{p}d\theta }<\infty .$
\end{defin}
The \textbf{Hardy space} $H^{\infty }({\mathbb{D}}^{n})$ is the
 space of holomorphic and bounded functions in ${\mathbb{D}}^{n}$
 equipped  with the sup norm.\ \par 
\quad The space $H^{p}({\mathbb{D}}^{n})$ possesses a \textbf{reproducing
 kernel} for any $a\in {\mathbb{D}}^{n},\ k_{a}(z)$:\ \par 
\quad \quad \quad $\displaystyle k_{a}(z)=\frac{1}{(1-\bar a_{1}z_{1})}{\times}\cdot
 \cdot \cdot {\times}\frac{1}{(1-\bar a_{n}z_{n})}.$\ \par 
We know that a function $f$ in $H^{p}$ has almost everywhere
 boundary values $f^{*}$ in $L^{p}({\mathbb{T}}^{n}).$ And we
 have the reproducing property:\ \par 
\quad \quad \quad $\displaystyle \forall f\in H^{p}({\mathbb{D}}^{n}),\ f(a):={\left\langle{f,k_{a}}\right\rangle}=\int_{{\mathbb{T}}^{n}}{f^{*}(e^{i\theta
 })\bar k_{a}(e^{i\theta })d\theta _{1}\cdot \cdot \cdot d\theta _{n}},$\ \par 
where ${\left\langle{\cdot ,\cdot }\right\rangle}$ is the scalar
 product of the Hilbert space $H^{2}({\mathbb{D}}^{n}).$\ \par 
\quad We shall use the notation:\ \par 
\quad \quad \quad $\displaystyle ((1-\bar az)):=(1-\bar a_{1}z_{1}){\times}\cdot
 \cdot \cdot {\times}(1-\bar a_{n}z_{n})$\ \par 
hence:\ \par 
\quad \quad \quad $\displaystyle ((1-\left\vert{a}\right\vert ^{2})):=(1-\left\vert{a_{1}}\right\vert
 ^{2}){\times}\cdot \cdot \cdot {\times}(1-\left\vert{a_{n}}\right\vert
 ^{2}).$\ \par 
Then the normalized reproducing kernel in $H^{p}({\mathbb{D}}^{n})$
 is, for $p\geq 1$:\ \par 
\quad \quad \quad $\displaystyle k_{a,p}(z)=\frac{((1-\left\vert{a}\right\vert
 ^{2}))^{1/p'}}{((1-\bar az))}.$\ \par 
\quad The same way for the unit ball ${\mathbb{B}}$ of ${\mathbb{C}}^{n},$
 we have:\ \par 

\begin{defin}
The {\bf Hardy space} $H^{p}({\mathbb{B}})$ is the set of holomorphic
 functions $f$ in ${\mathbb{B}}$ such that, with $z=r\zeta \in
 {\mathbb{B}},\ r\in (0,1),\ \zeta \in \partial {\mathbb{B}},$
 and $d\sigma $ the Lebesgue measure on $\partial {\mathbb{B}}$:\par 
\quad \quad \quad $\displaystyle {\left\Vert{f}\right\Vert}_{H^{p}}^{p}:=\sup _{r<1}\int_{\partial
 {\mathbb{B}}}{\left\vert{f(r\zeta )}\right\vert ^{p}d\sigma (\zeta )}<\infty .$
\end{defin}
The \textbf{Hardy space} $H^{\infty }({\mathbb{B}})$ is the space
 of holomorphic and bounded functions in ${\mathbb{B}}$ equipped
  with the sup norm.\ \par 
\quad The space $H^{p}({\mathbb{B}}),\ {\mathbb{B}}\subset {\mathbb{C}}^{n}$
 possesses a reproducing kernel for any $a\in {\mathbb{B}},\
 k_{a}(z)$ and a normalised one $k_{a,p}(z)$:\ \par 
\quad \quad \quad $\displaystyle k_{a}(z)=\frac{1}{(1-\bar a\cdot z)^{n}},\ k_{a,p}(z)=\frac{k_{a}(z)}{{\left\Vert{k_{a}}\right\Vert}_{p}},$
 with $\displaystyle {\left\Vert{k_{a}}\right\Vert}_{p}\simeq
 (1-\left\vert{a}\right\vert ^{2})^{-n/p'},$\ \par 
where we use the notation $\bar a\cdot z:=\sum_{j=1}^{n}{\bar
 a_{j}z_{j}}.$\ \par 
Again we know that a function $f$ in $H^{p}$ has almost everywhere
 boundary values $f^{*}$ in $L^{p}(\partial {\mathbb{B}}).$ And we have\ \par 
\quad \quad \quad $\displaystyle \forall f\in H^{p}({\mathbb{B}}),\ f(a):={\left\langle{f,k_{a}}\right\rangle}=\int_{\partial
 {\mathbb{B}}}{f(\zeta )\bar k_{a}(\zeta )d\zeta },$ \ \par 
where ${\left\langle{\cdot ,\cdot }\right\rangle}$ is the scalar
 product of the Hilbert space $H^{2}({\mathbb{B}}).$\ \par 
\ \par 
\quad Let $\Omega $ be either the ball ${\mathbb{B}}$ or the polydisc
 ${\mathbb{D}}^{n}.$ We shall use the notation: $\forall a\in
 \Omega ,\ \chi _{a}:={\left\Vert{k_{a}}\right\Vert}_{H^{2}}^{-2}.$
 Hence in ${\mathbb{D}}^{n},\ \chi _{a}=((1-\left\vert{a}\right\vert
 ^{2}))$ and in ${\mathbb{B}},\ \chi _{a}=(1-\left\vert{a}\right\vert
 ^{2})^{n}.$\ \par 

\begin{defin}
Let $S$ be a sequence of points in $\Omega .$ We say that $S$
 is a {\bf dual bounded sequence} in $H^{p}(\Omega )$ if there
 is a sequence $\lbrace \rho _{a}\rbrace _{a\in S}\subset H^{p}(\Omega
 )$ such that, with $\displaystyle k_{b,p'}$ the normalised reproducing
 kernel in $H^{p'}(\Omega )$ for the point $b\in \Omega $ and
 $\delta _{a,b}:=0$ for $a\neq b$ and $\displaystyle \delta _{a,b}:=1$
 for $a=b$:\par 
\quad \quad \quad $\displaystyle \exists \ C>0,\ \forall a,b\in S,\ {\left\langle{\rho
 _{a},k_{b,p'}}\right\rangle}=\delta _{a,b}$ and $\ {\left\Vert{\rho
 _{a}}\right\Vert}_{H^{p}(\Omega )}\leq C.$
\end{defin}
\quad In Section~\ref{e0} we shall consider the Grammian associated
 to a sequence of points $S$ in $\Omega .$ This is the infinite
 matrix $G$ given by\ \par 
\quad \quad \quad $\displaystyle G_{a,b}:={\left\langle{k_{a,p'},\ k_{b,p}}\right\rangle},$\ \par 
where the $\displaystyle k_{a,p}$ are the normalised reproducing
 kernels for $H^{p}(\Omega ).$\ \par 
\quad We consider $G$ as an operator on $???^{p'}(S)$ as\ \par 
\quad \quad \quad $\forall \mu \in ???^{p'}(S),\ G\mu :=\lbrace (G\mu )_{b}\rbrace
 _{b\in S},\ (G\mu )_{b}:=\sum_{a\in S}{\mu _{a}G_{a,b}}={\left\langle{\sum_{a\in
 S}{\mu _{a}k_{a,p'},k_{b,p}}}\right\rangle}.$\ \par 
In order to state some results, we shall need:\ \par 

\begin{defin}
Let $S$ be a sequence of points in $\Omega .$ We shall say that
 $S$ is {\bf strictly }$H^{p}\ ${\bf  interpolating} if the Grammian
 $G$ is bounded below on $???^{p'}(S).$
\end{defin}
\quad In Section~\ref{lE20} we study the case of the weighted Bergman
 spaces of the ball with the following definition.\ \par 

\begin{defin}
Let $f$ be a holomorphic function in ${\mathbb{B}}_{n}\subset
 {\mathbb{C}}^{n}$ and $k\in {\mathbb{N}};$ we say that $f\in
 A_{k}^{p}({\mathbb{B}})$ if\par 
\quad \quad \quad $\displaystyle \ {\left\Vert{f}\right\Vert}_{k,p}^{p}:=\int_{{\mathbb{B}}}{\left\vert{f(z)}\right\vert
 ^{p}(1-\left\vert{z}\right\vert ^{2})^{k}\,dm(z)}<\infty .$
\end{defin}
Where $\,dm$ is the Lebesgue measure in ${\mathbb{C}}^{n}.$\ \par 

\section{The results.}
\quad Recall that $\Omega $ denote either the ball ${\mathbb{B}}$ or
 the polydisc ${\mathbb{D}}^{n}.$\ \par 

\begin{thm}
Let $S$ be a sequence of points in $\Omega .$ Suppose that $S$
 is $H^{p}(\Omega )$ interpolating with $p\geq 2.$ Then $S$ has the LEP.
\end{thm}
Then we get:\ \par 

\begin{thm}
Let $S$ be a sequence of points in $\Omega .$ Suppose that $S$
 is Carleson and $S$ is strictly $H^{p}$ interpolating.\par 
Then $S$ is $H^{p}$ interpolating and has the LEP with range
 $E_{S}^{p}:=\mathrm{S}\mathrm{p}\mathrm{a}\mathrm{n}\lbrace
 k_{a,p},\ a\in S\rbrace .$
\end{thm}
And its converse:\ \par 

\begin{thm}
Let $S$ be a sequence of points in $\Omega .$ Suppose that $S$
 is $H^{p}$ interpolating and has the LEP with range $E_{S}^{p}:=\mathrm{S}\mathrm{p}\mathrm{a}\mathrm{n}\lbrace
 k_{a,p},\ a\in S\rbrace $ for a $p>2.$\par 
Then $S$ is strictly $H^{p}$ interpolating.
\end{thm}
\quad In Section~\ref{lE20} we study the case of the weighted Bergman
 spaces of the ball. By applying a "Subordination Lemma" in~\cite{amBerg78}
 Section 2, p.716\textbf{ }or, for a general form of it~\cite{subPrinAmar12},
 we get the same results as for the Hardy spaces.\ \par 

\begin{thm}
Let $S$ be a sequence of points in ${\mathbb{B}}\subset {\mathbb{C}}^{n}.$
 Suppose that $S$ is $A_{k}^{p}$-interpolating in ${\mathbb{B}}$
 for a $p\geq 2.$ then $S$ has the LEP with range $E_{S}^{p}:=\mathrm{S}\mathrm{p}\mathrm{a}\mathrm{n}\lbrace
 k_{a,p},\ a\in S\rbrace .$
\end{thm}

\begin{thm}
Let $S$ be a sequence of points in ${\mathbb{B}}.$ Suppose that
 $S$ is Carleson for $A_{k}^{p},$ and $S$ is strictly $A_{k}^{p}$
 interpolating.\par 
Then $S$ is $A_{k}^{p}$ interpolating and has the LEP with range
 $E_{S}^{p}:=\mathrm{S}\mathrm{p}\mathrm{a}\mathrm{n}\lbrace
 k_{a,p},\ a\in S\rbrace .$
\end{thm}
And its converse:\ \par 

\begin{thm}
Let $S$ be a sequence of points in ${\mathbb{B}}.$ Suppose that
 $S$ is $A_{k}^{p}$ interpolating and has the LEP with range
 $E_{S}^{p}:=\mathrm{S}\mathrm{p}\mathrm{a}\mathrm{n}\lbrace
 k_{a,p},\ a\in S\rbrace $ for a $p>2.$\par 
Then $S$ is strictly $A_{k}^{p}$ interpolating.
\end{thm}
\quad Finally we prove Babenko examples made of holomorphic functions.
 With $\Omega ={\mathbb{B}}_{2}$ or $\Omega ={\mathbb{D}}^{2}$ we get:\ \par 

\begin{thm}
There are sequences of reproducing kernels in $H^{p}(\Omega )$
 which are $p$-hilbertian but not  $p$-Riesz sequences.
\end{thm}

\section{Equivalence with the functional definitions.}

\begin{lem}
~\label{pH1} Let $S$ be a sequence of points in $\Omega .$ Then
 $S$ is $H^{p}(\Omega )$ interpolating iff:\par 
\quad \quad \quad $\displaystyle R_{p}H^{p}\supset ???^{p}(S).$
\end{lem}
\quad Proof.\ \par 
$\bullet $ First it is well known that $\displaystyle R_{p}H^{p}\supset
 ???^{p}(S)$ implies:\ \par 
\quad \quad \quad \begin{equation}  \exists c>0,\ \forall f\in \ H^{p}(\Omega ),\
 {\left\Vert{R_{p}f}\right\Vert}_{???^{p}(S)}\geq c{\left\Vert{f}\right\Vert}_{H^{p}(\Omega
 )}.\label{eI0}\end{equation}\ \par 
To see this consider the subspace:\ \par 
\quad \quad \quad $L_{K}:=\lbrace \lambda \in ???^{p}(S),\ \exists f\in H^{p},\
 {\left\Vert{f}\right\Vert}_{p}\leq K{\left\Vert{\lambda }\right\Vert}_{p},\
 R_{p}f=\lambda \rbrace .$\ \par 
Clearly $L_{K}$ is a closed subspace of $\displaystyle ???^{p}(S)$
 and we have $\bigcup_{K\in {\mathbb{N}}}{L_{K}}=???^{p}(S).$
 Hence by Baire's Theorem one of the $L_{K}$ contains a non void
 ball, which gives ~(\ref{eI0}).\ \par 
\ \par 
$\bullet $ Second suppose that~(\ref{eI0}) is true and $R_{p}H^{p}\supset
 ???^{p}(S).$ Then:\ \par 
\quad \quad \quad $\displaystyle \forall f\in H^{p},\ {\left\Vert{f}\right\Vert}_{p}\leq
 1,\ {\left\Vert{\sum_{a\in S}{\mu _{a}k_{a,p'}}}\right\Vert}_{p'}\geq
 \left\vert{{\left\langle{f,\sum_{a\in S}{\mu _{a}k_{a,p'}}}\right\rangle}}\right\vert
 .$\ \par 
But\ \par 
\quad \quad \quad $\displaystyle {\left\langle{f,\sum_{a\in S}{\mu _{a}k_{a,p'}}}\right\rangle}=\sum_{a\in
 S}{\bar \mu _{a}f(a)\chi _{a}^{1/p}}.$\ \par 
We take $f\in H^{p}$ such that  $f(a)\chi _{a,p'}=\lambda _{a}\in
 ???^{p}(S)$ with $\sum_{a\in S}{\bar \mu _{a}f(a)\chi _{a,p'}}={\left\Vert{\mu
 }\right\Vert}_{???^{p'}}.$ This is possible because $\displaystyle
 R_{p}H^{p}\supset ???^{p}(S)$ with ${\left\Vert{f}\right\Vert}_{p}\leq
 C{\left\Vert{\lambda }\right\Vert}_{???^{p}}.$ So we get, with
 ${\left\Vert{\lambda }\right\Vert}_{???^{p}}=1,$\ \par 
\quad \quad \quad $\displaystyle {\left\Vert{\sum_{a\in S}{\mu _{a}k_{a,p'}}}\right\Vert}_{p'}\geq
 \frac{1}{C}{\left\Vert{\mu }\right\Vert}_{???^{p'}}.$\ \par 
Which means that $\lbrace k_{a,p'}\rbrace _{a\in S}$ is $p$-besselian,
 hence $S$ is $H^{p}(\Omega )$ interpolating.\ \par 
\ \par 
$\bullet $ Third suppose that $S$ is $H^{p}(\Omega )$ interpolating,
 i.e. ${\left\Vert{\sum_{a\in S}{\mu _{a}k_{a,p'}}}\right\Vert}_{p'}\gtrsim
 {\left\Vert{\mu }\right\Vert}_{???^{p'}}.$\ \par 
We truncate $S$ to its first $N$ terms, call it $S_{N}.$ Then
 we have that $\lbrace k_{a,p'}\rbrace _{a\in S_{N}}$ has a $H^{p}$
 dual  sequence $\lbrace \rho _{a}\rbrace _{a\in S_{N}}\subset
 E_{S}^{p}:=\mathrm{S}\mathrm{p}\mathrm{a}\mathrm{n}\lbrace k_{a,p},\
 a\in S_{N}\rbrace ,$ i.e.\ \par 
\quad \quad \quad $\forall a\in S_{N},\ \exists \rho _{a}\in H^{p}::{\left\langle{\rho
 _{a},k_{b,p'}}\right\rangle}=\delta _{ab}.$\ \par 
\ \par 
Now we have that the dual space of $E_{S_{N}}^{p'}:=\mathrm{S}\mathrm{p}\mathrm{a}\mathrm{n}\lbrace
 k_{a,p'},\ a\in S_{N}\rbrace $ is the quotient space $H^{p}/\mathrm{A}\mathrm{n}\mathrm{n}(E_{S_{N}}^{p'}).$
 By definition the annihilator of $E_{S_{N}}^{p'}$ is $Z_{S_{N}}^{p}:=\lbrace
 g\in H^{p}::\forall a\in S_{N},\ g(a)=0\rbrace .$ So we have\ \par 
\quad \quad \quad $\displaystyle \forall \lambda \in ???^{p}(S_{N}),\ {\left\Vert{\sum_{a\in
 S_{N}}{\lambda _{a}\rho _{a}}}\right\Vert}_{H^{p}/Z_{S_{N}}^{p}}=\sup
 _{f\in E_{S_{N}}^{p'},{\left\Vert{f}\right\Vert}_{p'}\leq 1}\left\vert{{\left\langle{\sum_{a\in
 S_{N}}{\lambda _{a}\rho _{a}},f}\right\rangle}}\right\vert .$\ \par 
But $f\in E_{S}^{p'}$ means that $f=\sum_{b\in S_{N}}{\mu _{b}k_{b,p'}}$
 hence\ \par 
\quad \quad \quad $\displaystyle {\left\Vert{\sum_{a\in S_{N}}{\lambda _{a}\rho
 _{a}}}\right\Vert}_{H^{p}/Z_{S_{N}}^{p}}=\sup _{f\in E_{S_{N}}^{p'},{\left\Vert{f}\right\Vert}_{p'}\leq
 1}\left\vert{{\left\langle{\sum_{a\in S_{N}}{\lambda _{a}\rho
 _{a}},\sum_{b\in S_{N}}{\mu _{b}k_{b,p'}}}\right\rangle}}\right\vert .$\ \par 
So, using H\"older inequalities,\ \par 
\quad \quad \quad $\displaystyle {\left\Vert{\sum_{a\in S_{N}}{\lambda _{a}\rho
 _{a}}}\right\Vert}_{H^{p}/Z_{S_{N}}^{p}}=\sup _{f\in E_{S_{N}}^{p'},{\left\Vert{f}\right\Vert}_{p'}\leq
 1}\left\vert{\sum_{a\in S_{N}}{\lambda _{a}\bar \mu _{a}}}\right\vert
 \leq {\left\Vert{\lambda }\right\Vert}_{???^{p}}{\left\Vert{\mu
 }\right\Vert}_{???^{p'}}$\ \par 
by assumption, ${\left\Vert{\mu }\right\Vert}_{???^{p'}}\leq
 C{\left\Vert{\sum_{a\in S_{N}}{\mu _{a}k_{a,p'}}}\right\Vert}_{p'}=C{\left\Vert{f}\right\Vert}_{p'},$
 with $C$ independent of $N,$ hence\ \par 
\quad \quad \quad $\displaystyle {\left\Vert{\sum_{a\in S_{N}}{\lambda _{a}\rho
 _{a}}}\right\Vert}_{H^{p}/Z_{S_{N}}^{p}}\leq C{\left\Vert{\lambda
 }\right\Vert}_{???^{p}}.$\ \par 
Now the functions $g:=\sum_{a\in S_{N}}{\lambda _{a}\rho _{a}}+Z_{S}^{p}$
 interpolate the right values: $\forall a\in S,\ g(a)=\lambda
 _{a}\chi _{a,p'}.$ So we can find one $f$ in $H^{p}$ with norm
 less than $(C+\epsilon ){\left\Vert{\lambda }\right\Vert}_{???^{p}(S_{N})}.$\
 \par 
We proved that for any $N$:\ \par 
\quad \quad \quad $\displaystyle \forall \lambda \in ???^{p}(S_{N}),\ \exists f\in
 H^{p},\ {\left\Vert{f}\right\Vert}_{p}\leq (C+1){\left\Vert{\lambda
 }\right\Vert}_{???^{p}(S_{N})}$ with $\forall a\in S,\ {\left\langle{f,k_{a,p'}}\right\rangle}=\lambda
 _{a}.$\ \par 
Because $\displaystyle ???^{p}(S_{N})$ is dense in $\displaystyle
 ???^{p}(S)$ and none of the constants depend on $N,$ we get
 the result by letting $N\rightarrow \infty .$ I.e. the restriction
 operator verifies~(\ref{eI0}). $\blacksquare $\ \par 

\begin{lem}
~\label{eI1}Let $S$ be a sequence of points in $\Omega .$ Then
 $S$ is Carleson iff there exists a $p\in (1,\infty )$ such that:\par 
\quad \quad \quad $\displaystyle (*)\ \ \ \ \ \ C>0,\ \forall f\in \ H^{p}(\Omega
 ),\ {\left\Vert{R_{p}f}\right\Vert}_{???^{p}(S)}\leq C{\left\Vert{f}\right\Vert}_{H^{p}(\Omega
 )},$\par 
i.e., with the measure $d\mu :=\sum_{a\in S}{\chi _{a,p'}\delta _{a}},$\par 
\quad \quad \quad $\forall f\in H^{p},\ \int_{\Omega }{\left\vert{f}\right\vert
 ^{p}d\mu }\leq C^{p}{\left\Vert{f}\right\Vert}_{p}^{p}.$ \par 
Which means that $\mu $ is a Carleson measure in $\Omega .$
\end{lem}
\quad Proof.\ \par 
$\bullet $ Suppose first that $S$ verifies $(*)$ and take any
 $p\in (1,\infty ).$\ \par 
We have, for any $\mu \in ???^{p'}(S),$\ \par 
\quad \quad \quad $\displaystyle {\left\Vert{\sum_{a\in S}{\mu _{a}k_{a,p'}}}\right\Vert}_{p'}=\sup
 _{f\in H^{p},{\left\Vert{f}\right\Vert}_{p}\leq 1}\left\vert{{\left\langle{f,\sum_{a\in
 S}{\mu _{a}k_{a,p'}}}\right\rangle}}\right\vert .$\ \par 
But  $\left\vert{{\left\langle{f,\sum_{a\in S}{\mu _{a}k_{a,p'}}}\right\rangle}}\right\vert
 =\left\vert{\sum_{a\in S}{\bar \mu _{a}{\left\langle{f,k_{a,p'}}\right\rangle}}}\right\vert
 $ and by H\"older inequalities we get:\ \par 
\quad \quad \quad \begin{equation} \left\vert{ \sum_{a\in S}{\bar \mu _{a}{\left\langle{f,k_{a,p'}}\right\rangle}}}\right\vert
 \leq {\left({\sum_{a\in S}{\left\vert{{\left\langle{f,k_{a,p'}}\right\rangle}}\right\vert
 ^{p}}}\right)}^{1/p}{\left\Vert{\mu }\right\Vert}_{???^{p'}(S)}.\label{lE3}\end{equation}\
 \par 
$S$ Carleson means that:\ \par 
\quad \quad \quad $\displaystyle \forall f\in H^{p},\ \sum_{a\in S}{\left\vert{{\left\langle{f,k_{a,p'}}\right\rangle}}\right\vert
 ^{p}}\lesssim {\left\Vert{f}\right\Vert}_{p}^{p},$\ \par 
so we get by~(\ref{lE3}) that\ \par 
\quad \quad \quad $\displaystyle \left\vert{\sum_{a\in S}{\bar \mu _{a}{\left\langle{f,k_{a,p'}}\right\rangle}}}\right\vert
 \lesssim {\left\Vert{f}\right\Vert}_{p}^{p}{\left\Vert{\mu }\right\Vert}_{???^{p'}(S)}$\
 \par 
which gives ${\left\Vert{\sum_{a\in S}{\mu _{a}k_{a,p'}}}\right\Vert}_{p'}\lesssim
 {\left\Vert{\mu }\right\Vert}_{???^{p'}}$ because ${\left\Vert{f}\right\Vert}_{p}\leq
 1.$\ \par 
\ \par 
$\bullet $ Suppose now that\ \par 
\quad \quad \quad $\displaystyle {\left\Vert{\sum_{a\in S}{\mu _{a}k_{a,p'}}}\right\Vert}_{p'}\lesssim
 {\left\Vert{\mu }\right\Vert}_{???^{p'}}.$\ \par 
Take $f\in H^{p}$ then\ \par 
\quad \quad \quad $\displaystyle \left\vert{{\left\langle{f,\sum_{a\in S}{\mu _{a}k_{a,p'}}}\right\rangle}}\right\vert
 =\left\vert{\sum_{a\in S}{\bar \mu _{a}{\left\langle{f,k_{a,p'}}\right\rangle}}}\right\vert
 .$\ \par 
Now take $\mu \in ???^{p'}(S),\ {\left\Vert{\mu }\right\Vert}_{???^{p'}(S)}=1$
 such that $\left\vert{\sum_{a\in S}{\bar \mu _{a}{\left\langle{f,k_{a,p'}}\right\rangle}}}\right\vert
 ={\left({\sum_{a\in S}{\left\vert{{\left\langle{f,k_{a,p'}}\right\rangle}}\right\vert
 ^{p}}}\right)}^{1/p}.$ Then we get\ \par 
\quad \quad \quad $\displaystyle {\left({\sum_{a\in S}{\left\vert{{\left\langle{f,k_{a,p'}}\right\rangle}}\right\vert
 ^{p}}}\right)}^{1/p}=\left\vert{{\left\langle{f,\sum_{a\in S}{\mu
 _{a}k_{a,p'}}}\right\rangle}}\right\vert \leq {\left\Vert{f}\right\Vert}_{p}{\left\Vert{\sum_{a\in
 S}{\mu _{a}k_{a,p'}}}\right\Vert}_{p'}\lesssim {\left\Vert{f}\right\Vert}_{p}{\left\Vert{\mu
 }\right\Vert}_{???^{p'}}\lesssim {\left\Vert{f}\right\Vert}_{p}$\ \par 
because ${\left\Vert{\mu }\right\Vert}_{???^{p'}(S)}=1.$ This
 means that $S$ is verifies $(*).$ $\blacksquare $\ \par 

\begin{rem}
Recall that if $S$ is Carleson for a $p\in \lbrack 1,\infty \lbrack
 ,$ then it is Carleson for any $q\in \lbrack 1,\infty \lbrack
 $ as the geometric characterisation of Carleson measure, done
 for ${\mathbb{B}}$ by Hastings ~\cite{Hastings75} and done for
 ${\mathbb{D}}^{n}$ by Chang~\cite{AChang79}, proved.
\end{rem}

\section{Proof of the main result.}
\quad $S$ has the LEP for $H^{p}$ means that there is a bounded operator
 $E:???^{p}(S)\rightarrow H^{p}$ such that $\forall a\in S,\
 {\left\langle{E(\lambda ),\ k_{a,p'}}\right\rangle}=\lambda
 _{a}.$ We shall say that $S$ has the LEP \emph{with range} $F:=\mathrm{R}\mathrm{a}\mathrm{n}\mathrm{g}\mathrm{e}(E).$\
 \par 

\begin{thm}
~\label{lE17}Let $S$ be a sequence of points in $\Omega .$ Suppose
 that $S$ is $H^{p}$-interpolating in $\Omega $ for a $p\geq
 2.$ then $S$ has the LEP with range $E_{S}^{p}:=\mathrm{S}\mathrm{p}\mathrm{a}\mathrm{n}\lbrace
 k_{a,p},\ a\in S\rbrace .$
\end{thm}
\quad Proof.\ \par 
$\bullet $ The case $p=2$ is easy.\ \par 
Set $E_{S}^{2}:=\mathrm{S}\mathrm{p}\mathrm{a}\mathrm{n}\lbrace
 k_{a,2},\ a\in S\rbrace $ and $P$ the orthogonal projection
 from $H^{2}$ onto $E_{S}^{2}.$ Let $\lambda \in ???^{2}(S)$
 then, because $S$ is $H^{2}$ interpolating, there is a function
 $f\in H^{2}$ such that $\forall a\in S,\ {\left\langle{f,k_{a,2}}\right\rangle}=\lambda
 _{a}$ and ${\left\Vert{f}\right\Vert}_{2}\lesssim {\left\Vert{\lambda
 }\right\Vert}_{???^{2}(S)}.$ Set $g:=Pf$ then\ \par 
\quad \quad \quad $\displaystyle \forall a\in S,\ {\left\langle{g,k_{a,2}}\right\rangle}={\left\langle{Pf,k_{a,2}}\right\rangle}={\left\langle{f,Pk_{a,2}}\right\rangle}={\left\langle{f,k_{a,2}}\right\rangle}=\lambda
 _{a}$\ \par 
and ${\left\Vert{g}\right\Vert}_{2}={\left\Vert{Pf}\right\Vert}_{2}\leq
 {\left\Vert{f}\right\Vert}_{2}\lesssim {\left\Vert{\lambda }\right\Vert}_{???^{2}(S)}.$
 If there is another function $h\in E_{S}^{2}$ such that $\forall
 a\in S,\ {\left\langle{h,k_{a,2}}\right\rangle}=\lambda _{a}$
 then ${\left\langle{g-h,k_{a,2}}\right\rangle}=0$ hence $g-h\perp
 E_{S}^{2}\Rightarrow g-h=0.$ So the operator $E:\ \lambda \in
 ???^{2}(S)\rightarrow g\in E_{S}^{2}$ is the extension we are
 searching for.\ \par 
\ \par 
$\bullet $ The case $p>2.$\ \par 
Using the Theorem~\ref{lE7} we have that, if $\Omega ={\mathbb{B}},\
 S$ is Carleson. By Theorem 1.11 p. 3 in~\cite{CarlesonInter19}
 we have that, if $\Omega ={\mathbb{D}}^{n},\ S$ is Carleson.
 So in any cases we have that $S$ is Carleson.\ \par 
Because $S$ is $H^{p}$ interpolating, we have:\ \par 
\quad \quad \quad $\displaystyle \forall \mu \in ???^{p'}(S),\ {\left\Vert{\sum_{a\in
 S}{\mu _{a}k_{a,p'}}}\right\Vert}_{p'}\gtrsim {\left\Vert{\mu
 }\right\Vert}_{???^{p'}}.$\ \par 
Because $S$ is Carleson, we have:\ \par 
\quad \quad \quad $\displaystyle \forall \mu \in ???^{p'}(S),\ {\left\Vert{\sum_{a\in
 S}{\mu _{a}k_{a,p'}}}\right\Vert}_{p'}\lesssim {\left\Vert{\mu
 }\right\Vert}_{???^{p'}}.$\ \par 
Hence we get that the sequence $\lbrace k_{a,p'}\rbrace _{a\in
 S}$ is a $p$-Riesz basis for the space $E_{S}^{p'}:=\mathrm{S}\mathrm{p}\mathrm{a}\mathrm{n}\lbrace
 k_{a,p'},\ a\in S\rbrace .$\ \par 
This easily implies, see for instance Theorem 4.6, p. 15 in~\cite{ChalChev18},
 that there is a bounded operator $Q_{p'}:E_{S}^{p'}\rightarrow
 ???^{p'}(S)$ with bounded inverse such that $\forall a\in S,\
 \epsilon _{a}=Q_{p'}k_{a,p'},$ where $\epsilon _{a}$ is the
 canonical basis of $???^{p'}(S).$\ \par 
\quad Using Theorem 1.5 and Theorem 1.6 p. 179 in~\cite{AmarIntInt07},
 we get that if $S$ is $H^{p}$ interpolating then $S$ is $H^{q}$
 interpolating for $q<p.$ Because $p>2$ we get that $p'<2<p$
 hence $S$ is also $H^{p'}$ interpolating.\ \par 
\quad So we can apply the above inequalities with $p'$ instead of $p$
 to get that the sequence $\lbrace k_{a,p}\rbrace _{a\in S}$
 is a $p$-Riesz basis for the space $E_{S}^{p}:=\mathrm{S}\mathrm{p}\mathrm{a}\mathrm{n}\lbrace
 k_{a,p}\ a\in S\rbrace .$ And again there is an bounded operator
 $Q_{p}:E_{S}^{p}\rightarrow ???^{p}(S)$ with bounded inverse
 such that $\forall a\in S,\ f_{a}=Q_{p}k_{a,p}$ where $f_{a}$
 is the canonical basis of $???^{p}(S).$\ \par 
\quad The idea is to extend the Proposition 2 p. 13 in~\cite{AmarThesis77}
 done for Hilbert spaces to our case. We get\ \par 
\quad \quad \quad $\displaystyle \delta _{ab}={\left\langle{\epsilon _{a},f_{b}}\right\rangle}={\left\langle{Q_{p'}k_{a,p'},\
 Q_{p}k_{a,p}}\right\rangle}={\left\langle{k_{a,p'},\ (Q_{p'})^{\ast
 }Q_{p}k_{a,p}}\right\rangle}.$\ \par 
Hence, setting $\forall a\in S,\ \rho _{a}:=(Q_{p'})^{\ast }Q_{p}k_{a,p},$
 we get that $\lbrace \rho _{a}\rbrace _{a\in S}$ is a bounded
 dual sequence to $\lbrace k_{a,p'}\rbrace _{a\in S}$ contained
 in $E_{S}^{p}.$ I.e.\ \par 
\quad \quad \quad $\displaystyle \forall a\in S,\ {\left\langle{\rho _{a},k_{b,p'}}\right\rangle}=\delta
 _{ab}$ and $\exists C>0,\ \forall a\in S,\ {\left\Vert{\rho
 _{a}}\right\Vert}_{p}\leq C.$\ \par 
Moreover we have:\ \par 
\quad \quad \quad $\displaystyle \forall \lambda \in ???^{p}(S),\ {\left\Vert{\sum_{a\in
 S}{\lambda _{a}\rho _{a}}}\right\Vert}_{p}={\left\Vert{(Q_{p'})^{\ast
 }Q_{p}\sum_{a\in S}{\lambda _{a}k_{a,p}}}\right\Vert}_{p}.$\ \par 
But, using that $S$ is Carleson, we get\ \par 
\quad \quad \quad $\displaystyle {\left\Vert{\sum_{a\in S}{\lambda _{a}k_{a,p}}}\right\Vert}_{p}\lesssim
 {\left\Vert{\lambda }\right\Vert}_{???^{p}(S)}.$\ \par 
So\ \par 
\quad \quad \quad $\displaystyle \forall \lambda \in ???^{p}(S),\ {\left\Vert{\sum_{a\in
 S}{\lambda _{a}\rho _{a}}}\right\Vert}_{p}={\left\Vert{(Q_{p'})^{\ast
 }Q_{p}\sum_{a\in S}{\lambda _{a}k_{a,p}}}\right\Vert}_{p}\leq $\ \par 
\quad \quad \quad \quad \quad \quad \quad \quad \quad $\displaystyle \leq {\left\Vert{Q_{p'}}\right\Vert}{\left\Vert{Q_{p}}\right\Vert}{\left\Vert{\sum_{a\in
 S}{\lambda _{a}k_{a,p}}}\right\Vert}_{p}\lesssim {\left\Vert{Q_{p'}}\right\Vert}{\left\Vert{Q_{p}}\right\Vert}{\left\Vert{\lambda
 }\right\Vert}_{p}.$\ \par 
Because $Q_{p}$ and $Q_{p'}$ are bounded. We finally get:\ \par 
\quad \quad \quad $\displaystyle \forall \lambda \in ???^{p}(S),\ {\left\Vert{\sum_{a\in
 S}{\lambda _{a}\rho _{a}}}\right\Vert}_{p}\lesssim {\left\Vert{\lambda
 }\right\Vert}_{p}.$\ \par 
It remains to define the extension operator:\ \par 
\quad \quad \quad $\displaystyle \forall \lambda \in ???^{p}(S),\ E(\lambda )(z):=\sum_{a\in
 S}{\lambda _{a}\rho _{a}(z)}$\ \par 
to end the proof of the theorem. $\blacksquare $\ \par 

\begin{rem}
The Proposition 4.5 in~\cite{ChalChev18} says, in a fairly general
 situation, that $S$ is $H^{p}$ interpolating with the LEP iff
 we have that the sequence $\lbrace \rho _{a}\rbrace _{a\in S}$
 is $p$-hilbertian, i.e. iff ${\left\Vert{\sum_{a\in S}{\lambda
 _{a}\rho _{a}(z)}}\right\Vert}\leq C{\left\Vert{\lambda }\right\Vert}_{???^{p}(S)}.$
 This is pretty clear in our case here.
\end{rem}

\section{Strictly $H^{p}$ interpolating sequences.~\label{e0}}
\quad Let $S$ be a sequence of points in $\Omega .$ We consider the
 Grammian associated to $S.$ This is the infinite matrix $G$ given by\ \par 
\quad \quad \quad $\displaystyle G_{a,b}:={\left\langle{k_{a,p'},\ k_{b,p}}\right\rangle}.$\ \par 
We consider $G$ as an operator on $???^{p'}(S)$ defined as\ \par 
\quad \quad \quad $\displaystyle \forall \mu \in ???^{p'}(S),\ G\mu :=\lbrace (G\mu
 )_{b}\rbrace _{b\in S},\ (G\mu )_{b}:=\sum_{a\in S}{\mu _{a}G_{a,b}}={\left\langle{\sum_{a\in
 S}{\mu _{a}k_{a,p'},k_{b,p}}}\right\rangle}.$\ \par 
The adjoint $G^{\ast }$ of $G$ is the matrix $G^{\ast }_{a,b}=\bar
 G_{b,a}$ hence:\ \par 
\quad \quad \quad $\displaystyle G^{\ast }_{a,b}={\left\langle{k_{a,p},k_{b,p'}}\right\rangle}.$\
 \par 

\begin{rem}
~\label{lE11} We have the well known results:\par 
\quad \quad \quad $\displaystyle {\left\Vert{G}\right\Vert}_{???^{p}\rightarrow
 ???^{p}}={\left\Vert{G^{\ast }}\right\Vert}_{???^{p'}\rightarrow
 ???^{p'}}$\par 
and $G$ is bounded below, i.e. ${\left\Vert{G\mu }\right\Vert}_{???^{p'}}\gtrsim
 {\left\Vert{\mu }\right\Vert}_{???^{p'}},$ iff $G^{\ast }$ is bounded below.
\end{rem}
First we generalise to our setting half of the Proposition 9.5,
 p. 127 of~\cite{AglMCar02}:\ \par 

\begin{prop}
~\label{lE18}Let $S$ be a sequence of points in $\Omega .$ Let:\par 
\quad (CS) The sequence $S$ is Carleson.\par 
\quad (BG) The associated Grammian $G$ is bounded on $???^{p'}(S).$\par 
Then (CS) implies (BG).
\end{prop}
\quad Proof.\ \par 
We have\ \par 
\quad \quad \quad $\displaystyle (G\mu )_{b}:=\sum_{a\in S}{\mu _{a}{\left\langle{k_{a,p'},\
 k_{b,p}}\right\rangle}}={\left\langle{\sum_{a\in S}{\mu _{a}k_{a,p'}},\
 k_{b,p}}\right\rangle}.$\ \par 
Hence\ \par 
\quad \quad \quad $\displaystyle {\left\langle{G\mu ,\lambda }\right\rangle}=\sum_{b\in
 S}{\bar \lambda _{b}(G\mu )_{b}}={\left\langle{\sum_{a\in S}{\mu
 _{a}k_{a,p'}},\ \sum_{b\in S}{\lambda _{b}k_{b,p}}}\right\rangle}.$\ \par 
To get that $G$ is bounded on $???^{p'}(S)$ amounts to prove\ \par 
\quad \quad \quad $\displaystyle \left\vert{\sum_{a\in S}{\bar \lambda _{a}(G\mu
 )_{a}}}\right\vert \lesssim {\left\Vert{\lambda }\right\Vert}_{???^{p}(S)}{\left\Vert{\mu
 }\right\Vert}_{???^{p'}(S)}.$\ \par 
But, using H\"older inequalities,\ \par 
\quad \quad \quad \begin{equation} \left\vert{ {\left\langle{\sum_{a\in S}{\mu
 _{a}k_{a,p'}},\ \sum_{a\in S}{\lambda _{a}k_{a,p}}}\right\rangle}}\right\vert
 \leq {\left\Vert{\sum_{a\in S}{\mu _{a}k_{a,p'}}}\right\Vert}_{p'}{\left\Vert{\sum_{a\in
 S}{\lambda _{a}k_{a,p}}}\right\Vert}_{p}.\label{lE8}\end{equation}\ \par 
Suppose now that the sequence $S$ is Carleson. By definition this means:\ \par 
\quad \quad \quad $\displaystyle \forall \mu \in ???^{p'}(S),\ {\left\Vert{\sum_{a\in
 S}{\mu _{a}k_{a,p'}}}\right\Vert}_{p'}\lesssim {\left\Vert{\mu
 }\right\Vert}_{???^{p'}}$\ \par 
and\ \par 
\quad \quad \quad $\displaystyle \forall \lambda \in ???^{p}(S),\ {\left\Vert{\sum_{a\in
 S}{\lambda _{a}k_{a,p}}}\right\Vert}_{p'}\lesssim {\left\Vert{\lambda
 }\right\Vert}_{???^{p}}.$\ \par 
Hence~(\ref{lE8}) gives\ \par 
\quad \quad \quad $\displaystyle \left\vert{{\left\langle{G\mu ,\lambda }\right\rangle}}\right\vert
 =\left\vert{{\left\langle{\sum_{b\in S}{\mu _{b}k_{b,p'}},\
 \sum_{a\in S}{\lambda _{a}k_{a,p}}}\right\rangle}}\right\vert
 \lesssim {\left\Vert{\lambda }\right\Vert}_{???^{p}(S)}{\left\Vert{\mu
 }\right\Vert}_{???^{p'}(S)}.$\ \par 
Hence $G$ is bounded. The proof is complete. $\blacksquare $\ \par 

\begin{defin}
Let $S$ be a sequence of points in $\Omega .$ We shall say that
 $S$ is {\bf strictly }$H^{p}\ ${\bf  interpolating} if the operator
 $G$ is bounded below on $???^{p'}(S).$
\end{defin}

\begin{rem}
~\label{lE12} Using Remark~\ref{lE11}, we have that $S$ is strictly
 $H^{p}$ interpolating iff $S$ is strictly $H^{p'}$ interpolating.
\end{rem}

\begin{thm}
~\label{lE13}Let $S$ be a sequence of points in $\Omega .$ Suppose
 that $S$ is Carleson and $S$ is strictly $H^{p}$ interpolating.\par 
Then $S$ is $H^{p}$ interpolating and has the LEP with range
 $E_{S}^{p}:=\mathrm{S}\mathrm{p}\mathrm{a}\mathrm{n}\lbrace
 k_{a,p},\ a\in S\rbrace .$
\end{thm}
\quad For the proof we shall need a lemma, where $Z_{S}^{p'}:=\lbrace
 u\in H^{p'}::\forall a\in S,\ u(a)=0\rbrace .$\ \par 

\begin{lem}
~\label{lE9} If $S$ is strictly $H^{p}$ interpolating and Carleson
 then we have:\par 
\quad \quad \quad $\displaystyle \forall \mu \in ???^{p'}(S),\ \forall u\in Z_{S}^{p'},\
 \ {\left\Vert{u+\sum_{a\in S}{\mu _{a}k_{a,p'}}}\right\Vert}_{p'}\gtrsim
 {\left\Vert{\mu }\right\Vert}_{???^{p'}(S)}.$
\end{lem}
\quad Proof.\ \par 
We have, for any $f\in H^{p},$\ \par 
\quad \quad \quad $\displaystyle {\left\Vert{u+\sum_{a\in S}{\mu _{a}k_{a,p'}}}\right\Vert}_{p'}\geq
 \frac{\left\vert{{\left\langle{f,u}\right\rangle}+{\left\langle{f,\sum_{a\in
 S}{\mu _{a}k_{a,p'}}}\right\rangle}}\right\vert }{{\left\Vert{f}\right\Vert}_{p}}.$\
 \par 
We choose $f:=\sum_{a\in S}{\lambda _{a}k_{a,p}}$ then ${\left\langle{f,u}\right\rangle}=\sum_{a\in
 S}{\lambda _{a}{\left\langle{k_{a,p},u}\right\rangle}}=0$ because
 $\forall a\in S,\ u(a)=0.$\ \par 
Hence we have\ \par 
\quad \quad \quad $\displaystyle {\left\Vert{u+\sum_{a\in S}{\mu _{a}k_{a,p'}}}\right\Vert}_{p'}\geq
 \frac{\left\vert{{\left\langle{\sum_{a\in S}{\lambda _{a}k_{a,p}},\sum_{a\in
 S}{\mu _{a}k_{a,p'}}}\right\rangle}}\right\vert }{{\left\Vert{\sum_{a\in
 S}{\lambda _{a}k_{a,p}}}\right\Vert}_{p}}=\frac{\left\vert{{\left\langle{G\mu
 ,\lambda }\right\rangle}}\right\vert }{{\left\Vert{\sum_{a\in
 S}{\lambda _{a}k_{a,p}}}\right\Vert}_{p}}.$\ \par 
We choose $\lambda \in ???^{p}(S),\ {\left\Vert{\lambda }\right\Vert}_{???^{p}(S)}=1$
 such that ${\left\langle{G\mu ,\lambda }\right\rangle}={\left\Vert{G\mu
 }\right\Vert}_{???^{p'}(S)}.$ Now we use that $G$ is bounded
 below to get\ \par 
\quad \quad \quad $\displaystyle {\left\langle{G\mu ,\lambda }\right\rangle}={\left\Vert{G\mu
 }\right\Vert}_{???^{p'}(S)}\gtrsim {\left\Vert{\mu }\right\Vert}_{???^{p'}(S)}.$\
 \par 
Because $S$ is Carleson:\ \par 
\quad \quad \quad $\displaystyle \forall \lambda \in ???^{p}(S),\ {\left\Vert{\sum_{a\in
 S}{\lambda _{a}k_{a,p}}}\right\Vert}_{p}\lesssim {\left\Vert{\lambda
 }\right\Vert}_{???^{p}},$\ \par 
hence ${\left\Vert{\sum_{a\in S}{\lambda _{a}k_{a,p}}}\right\Vert}_{p}\lesssim
 {\left\Vert{\lambda }\right\Vert}_{???^{p}}=1.$\ \par 
\quad So we get\ \par 
\quad \quad \quad $\displaystyle {\left\Vert{u+\sum_{a\in S}{\mu _{a}k_{a,p'}}}\right\Vert}_{p'}\geq
 \frac{{\left\langle{G\mu ,\lambda }\right\rangle}}{{\left\Vert{\sum_{a\in
 S}{\lambda _{a}k_{a,p}}}\right\Vert}_{p}}\gtrsim {\left\Vert{\mu
 }\right\Vert}_{???^{p'}(S)}.$\ \par 
The proof of the lemma is complete. $\blacksquare $\ \par 

\begin{rem}
In fact this lemma proves that if $S$ is strictly $H^{p}$ interpolating
 and Carleson then we have:\par 
\quad \quad \quad $\displaystyle {\left\Vert{\sum_{a\in S}{\mu _{a}k_{a,p'}}}\right\Vert}_{H^{p'}/Z_{S}^{p'}}\gtrsim
 {\left\Vert{\mu }\right\Vert}_{???^{p'}(S)}.$
\end{rem}
\quad \textbf{Proof of the Theorem.}\ \par 
Suppose that $S$ is strictly $H^{p}$ interpolating.\ \par 
We truncate $S$ to its first $N$ terms, call it $S_{N}.$ Then
 we have that $\lbrace k_{a,p'}\rbrace _{a\in S_{N}}$ has a $H^{p}$
 dual bounded sequence $\lbrace \rho _{a}\rbrace _{a\in S_{N}}\subset
 E_{S}^{p}:=\mathrm{S}\mathrm{p}\mathrm{a}\mathrm{n}\lbrace k_{a,p},\
 a\in S_{N}\rbrace .$\ \par 
Take $\lambda \in ???^{p}(S_{N})$ we want to estimate  ${\left\Vert{\sum_{a\in
 S_{N}}{\lambda _{a}\rho _{a}}}\right\Vert}_{p}$ and, because
 the dual space of $E_{S_{N}}^{p}$ is $H^{p'}/Z_{S_{N}}^{p'}$ we have\ \par 
\quad \quad \quad $\displaystyle {\left\Vert{\sum_{a\in S_{N}}{\lambda _{a}\rho
 _{a}}}\right\Vert}_{p}=\sup _{f\in H^{p'}/Z_{S_{N}}^{p'}}\frac{\left\vert{{\left\langle{\sum_{a\in
 S_{N}}{\lambda _{a}\rho _{a}},\ f}\right\rangle}}\right\vert
 }{{\left\Vert{f}\right\Vert}_{p'}}.$\ \par 
\quad But $f\in H^{p'}/Z_{S_{N}}^{p'}$ can written, with $u\in Z_{S_{N}}^{p'},$\ \par 
\quad \quad \quad $\displaystyle f=u+\sum_{a\in S_{N}}{\mu _{a}k_{a,p'}}.$\ \par 
So we have to compute\ \par 
\quad \quad \quad $\displaystyle {\left\langle{\sum_{a\in S_{N}}{\lambda _{a}\rho
 _{a}},\ f}\right\rangle}={\left\langle{\sum_{a\in S_{N}}{\lambda
 _{a}\rho _{a}},\ u+\sum_{a\in S_{N}}{\mu _{a}k_{a,p'}}}\right\rangle}.$\ \par 
But $\sum_{a\in S_{N}}{\lambda _{a}\rho _{a}}\in E_{S}^{p}$ and
 $u\in Z_{S_{N}}^{p'},$ imply ${\left\langle{\sum_{a\in S_{N}}{\lambda
 _{a}\rho _{a}},\ u}\right\rangle}=0.$ So it remains\ \par 
\quad \quad \quad $\displaystyle {\left\langle{\sum_{a\in S_{N}}{\lambda _{a}\rho
 _{a}},\ \sum_{a\in S_{N}}{\mu _{a}k_{a,p'}}}\right\rangle}=\sum_{a\in
 S_{N}}{\lambda _{a}\bar \mu _{a}}$\ \par 
because ${\left\langle{\rho _{a},k_{b,p'}}\right\rangle}=\delta
 _{a,b}.$ So we get\ \par 
\quad \quad \quad $\displaystyle \left\vert{{\left\langle{\sum_{a\in S_{N}}{\lambda
 _{a}\rho _{a}},\ \sum_{a\in S_{N}}{\mu _{a}k_{a,p'}}}\right\rangle}}\right\vert
 =\left\vert{\sum_{a\in S_{N}}{\lambda _{a}\bar \mu _{a}}}\right\vert
 \leq {\left\Vert{\lambda }\right\Vert}_{???^{p}(S)}{\left\Vert{\mu
 }\right\Vert}_{???^{p'}(S)}.$\ \par 
Hence\ \par 
\quad \quad \quad \begin{equation} {\left\Vert{ \sum_{a\in S_{N}}{\lambda _{a}\rho
 _{a}}}\right\Vert}_{p}\leq \sup _{f\in H^{p'}/Z_{S_{N}}^{p'}}\frac{{\left\Vert{\lambda
 }\right\Vert}_{???^{p}(S_{N})}{\left\Vert{\mu }\right\Vert}_{???^{p'}(S_{N})}}{{\left\Vert{u+\sum_{a\in
 S_{N}}{\mu _{a}k_{a,p'}}}\right\Vert}_{p'}}.\label{lE10}\end{equation}\ \par 
Now we use Lemma~\ref{lE9} to get\ \par 
\quad \quad \quad $\displaystyle {\left\Vert{u+\sum_{a\in S_{N}}{\mu _{a}k_{a,p'}}}\right\Vert}_{p'}\gtrsim
 {\left\Vert{\mu }\right\Vert}_{???^{p'}(S_{N})}.$\ \par 
Putting it in~(\ref{lE10}) we deduce that\ \par 
\quad \quad \quad $\displaystyle {\left\Vert{\sum_{a\in S_{N}}{\lambda _{a}\rho
 _{a}}}\right\Vert}_{p}\lesssim {\left\Vert{\lambda }\right\Vert}_{???^{p}}.$\
 \par 
The constant under the $\lesssim $ is independent of $N,$ so
 we get, in particular, that ${\left\Vert{\rho _{a}^{N}}\right\Vert}\leq
 C.$ Now we use the diagonal process to let $N\rightarrow \infty
 $ and to get that there is a dual bounded sequence $\lbrace
 \rho _{a}\rbrace _{a\in S}\subset E_{S}^{p}$ such that\ \par 
\quad \quad \quad $\displaystyle \forall \lambda \in ???^{p}(S),\ {\left\Vert{\sum_{a\in
 S}{\lambda _{a}\rho _{a}}}\right\Vert}_{p}\lesssim {\left\Vert{\lambda
 }\right\Vert}_{???^{p}}.$\ \par 
It remains to define the extension operator:\ \par 
\quad \quad \quad $\displaystyle \forall \lambda \in ???^{p}(S),\ E(\lambda )(z):=\sum_{a\in
 S}{\lambda _{a}\rho _{a}(z)},$\ \par 
to end the proof of the theorem. $\blacksquare $\ \par 
\ \par 
\quad Using Remark~\ref{lE12} and Theorem~\ref{lE13}, we get that if
 $S$ is strictly $H^{p}$ interpolating, then $S$ is strictly
 $H^{p'}$ interpolating, hence $S$ is $H^{r}$ interpolating with
 the LEP for $r=\max (p,p')\geq 2.$\ \par 
\quad This leads to the converse of Theorem~\ref{lE13}.\ \par 

\begin{thm}
Let $S$ be a sequence of points in $\Omega .$ Suppose that $S$
 is $H^{p}$ interpolating and has the LEP with range $E_{S}^{p}:=\mathrm{S}\mathrm{p}\mathrm{a}\mathrm{n}\lbrace
 k_{a,p},\ a\in S\rbrace $ for a $p>2.$\par 
Then $S$ is strictly $H^{p}$ interpolating.
\end{thm}
\quad Proof.\ \par 
Because $p>2$ then $S$ is Carleson by Theorem~\ref{lE7} for the
 ball ${\mathbb{B}}$ and by Theorem 1.11 p. 3 in~\cite{CarlesonInter19}
 for the polydisc ${\mathbb{D}}^{n}.$\ \par 
\quad Using Theorem 1.5 and Theorem 1.6 p. 179 in~\cite{AmarIntInt07},
 we get that if $S$ is $H^{p}$ interpolating then $S$ is $H^{q}$
 interpolating for $q<p.$ Because $p>2$ we get that $p'<2<p$
 hence $S$ is also $H^{p'}$ interpolating.\ \par 
\ \par 
\quad Let ${\left\langle{\nu ,G\mu }\right\rangle}=\sum_{a,b\in S}{\nu
 _{a}\bar \mu _{b}{\left\langle{k_{a,p},\ k_{b,p'}}\right\rangle}}.$ \ \par 
$S$ being $H^{p}$ interpolating means ${\left\Vert{\sum_{a\in
 S}{\mu _{a}k_{a,p'}}}\right\Vert}_{p'}\gtrsim {\left\Vert{\mu
 }\right\Vert}_{???^{p'}}.$\ \par 
On the other hand, because $S$ is also $H^{p'}$ interpolating
 and Carleson, then:\ \par 
\quad \quad \quad \quad \quad \begin{equation} {\left\Vert{ \sum_{a\in S}{\nu _{a}k_{a,p}}}\right\Vert}_{p}\simeq
 {\left\Vert{\nu }\right\Vert}_{???^{p}}.\label{lE15}\end{equation}\ \par 
To get the norm of $\sum_{a\in S}{\mu _{a}k_{a,p'}}$ by duality,
 we have to test on $f\in H^{p}/Z_{S}^{p}.$ Such an $f$ can be written:\ \par 
\quad \quad \quad $\displaystyle f=u+\sum_{a\in S}{\nu _{a}k_{a,p}}$ with $u\in Z_{S}^{p}.$\ \par 
\ \par 
$S$ being Carleson means:\ \par 
\quad \quad \quad \begin{equation}  \forall f\in H^{p},\ {\left\Vert{\lbrace {\left\langle{f,k_{a,p'}}\right\rangle}\rbrace
 }\right\Vert}_{???^{p}(S)}\lesssim {\left\Vert{f}\right\Vert}_{p}.\label{lE14}\end{equation}\
 \par 
So, because $S$ is $H^{p}$ interpolating with the LEP with range
 $E_{S}^{p},$ we have that it exists a dual basis $\lbrace \rho
 _{a}\rbrace _{a\in S}\subset E_{S}^{p}$ such that:\ \par 
\quad \quad \quad $\displaystyle Pf:=\sum_{a\in S}{{\left\langle{f,k_{a,p'}}\right\rangle}\rho
 _{a}}\in E_{S}^{p}\subset H^{p}$ with $\displaystyle {\left\Vert{Pf}\right\Vert}_{p}\lesssim
 {\left\Vert{\lbrace {\left\langle{f,k_{a,p'}}\right\rangle}\rbrace
 }\right\Vert}_{???^{p}(S)}$\ \par 
hence  ${\left\Vert{Pf}\right\Vert}_{p}\lesssim {\left\Vert{f}\right\Vert}_{p}$
 by~(\ref{lE14}).\ \par 
Because $Pf=f$ on $S,$ we get also $P^{2}=P.$\ \par 
\quad Now\ \par 
\quad \quad \quad $\displaystyle {\left\langle{f,\sum_{a\in S}{\mu _{a}k_{a,p'}}}\right\rangle}={\left\langle{u+\sum_{a\in
 S}{\nu _{a}k_{a,p}},\sum_{a\in S}{\mu _{a}k_{a,p'}}}\right\rangle}={\left\langle{\sum_{a\in
 S}{\nu _{a}k_{a,p}},\sum_{a\in S}{\mu _{a}k_{a,p'}}}\right\rangle},$\ \par 
because ${\left\langle{u,\sum_{a\in S}{\mu _{a}k_{a,p'}}}\right\rangle}=0.$\
 \par 
Hence\ \par 
\quad \quad \quad $\displaystyle {\left\Vert{\sum_{a\in S}{\mu _{a}k_{a,p'}}}\right\Vert}_{p'}=\sup
 _{f\in H^{p}/Z_{S}^{p}}\frac{\left\vert{{\left\langle{\sum_{a\in
 S}{\nu _{a}k_{a,p}},\sum_{a\in S}{\mu _{a}k_{a,p'}}}\right\rangle}}\right\vert
 }{{\left\Vert{u+\sum_{a\in S}{\nu _{a}k_{a,p}}}\right\Vert}_{p}}.$\ \par 
But $P(u+\sum_{a\in S}{\nu _{a}k_{a,p}})=\sum_{a\in S}{\lambda
 _{a}\rho _{a}}$ with ${\left\Vert{\sum_{a\in S}{\lambda _{a}\rho
 _{a}}}\right\Vert}_{p}\lesssim {\left\Vert{u+\sum_{a\in S}{\nu
 _{a}k_{a,p}}}\right\Vert}_{p}\leq 1.$\ \par 
\quad Hence, because $\sum_{a\in S}{\nu _{a}k_{a,p}}=\sum_{a\in S}{\lambda
 _{a}\rho _{a}},$ we choose $f=u+\sum_{a\in S}{\nu _{a}k_{a,p}}$
 realizing the norm:\ \par 
\quad \quad \quad $\displaystyle {\left\Vert{\mu }\right\Vert}_{???^{p'}}\leq {\left\Vert{\sum_{a\in
 S}{\mu _{a}k_{a,p'}}}\right\Vert}_{p'}\lesssim \frac{\left\vert{{\left\langle{\sum_{a\in
 S}{\nu _{a}k_{a,p}},\sum_{a\in S}{\mu _{a}k_{a,p'}}}\right\rangle}}\right\vert
 }{{\left\Vert{\sum_{a\in S}{\lambda _{a}\rho _{a}}}\right\Vert}_{p}}.$\ \par 
But by~(\ref{lE15}) we get ${\left\Vert{\sum_{a\in S}{\nu _{a}k_{a,p}}}\right\Vert}_{p}\simeq
 {\left\Vert{\nu }\right\Vert}_{???^{p}}$ hence from $\sum_{a\in
 S}{\nu _{a}k_{a,p}}=\sum_{a\in S}{\lambda _{a}\rho _{a}}$ we get\ \par 
\quad \quad \quad $\displaystyle {\left\Vert{\sum_{a\in S}{\lambda _{a}\rho _{a}}}\right\Vert}_{p}={\left\Vert{\sum_{a\in
 S}{\nu _{a}k_{a,p}}}\right\Vert}_{p}\simeq {\left\Vert{\nu }\right\Vert}_{???^{p}}.$\
 \par 
Hence\ \par 
\quad \quad \quad $\displaystyle {\left\Vert{\mu }\right\Vert}_{???^{p'}}\leq {\left\Vert{\sum_{a\in
 S}{\mu _{a}k_{a,p'}}}\right\Vert}_{p'}\lesssim \frac{\left\vert{{\left\langle{\sum_{a\in
 S}{\nu _{a}k_{a,p}},\sum_{a\in S}{\mu _{a}k_{a,p'}}}\right\rangle}}\right\vert
 }{{\left\Vert{\nu }\right\Vert}_{???^{p}}}$\ \par 
which means exactly that $G$ is bounded below. Hence $S$ is strictly
 $H^{p}$ interpolating. $\blacksquare $\ \par 

\begin{rem}
In the case $p=2$ we retrieve the fact that $S$ is $H^{2}$ interpolating
 iff the Grammian $G$ is bounded below.
\end{rem}

\section{Bergman spaces of the ball.~\label{lE20}}
Let us recall the definition of the Bergman spaces of the ball
 we are interested in. We shall use the notation $\bar a\cdot
 z:=\sum_{j=1}^{n}{\bar a_{j}z_{j}}$ hence $\left\vert{z}\right\vert
 ^{2}:=z\cdot \bar z.$\ \par 

\begin{defin}
Let $f$ be a holomorphic function in ${\mathbb{B}}_{n}\subset
 {\mathbb{C}}^{n}$ and $k\in {\mathbb{N}};$ we say that $f\in
 A_{k}^{p}({\mathbb{B}})$ if\par 
\quad \quad \quad $\displaystyle \ {\left\Vert{f}\right\Vert}_{k,p}^{p}:=\int_{{\mathbb{B}}}{\left\vert{f(z)}\right\vert
 ^{p}(1-\left\vert{z}\right\vert ^{2})^{k}\,dm(z)}<\infty .$
\end{defin}
Where $\,dm$ is the Lebesgue measure in ${\mathbb{C}}^{n}.$\ \par 
\quad The space $A_{k}^{p}({\mathbb{B}})$ possesses a reproducing kernel
 for any $a\in {\mathbb{B}},\ k_{a}(z)$ and a normalised one
 $k_{a,p}(z)$:\ \par 
\quad \quad \quad $\displaystyle k_{a}(z)=\frac{1}{(1-\bar a\cdot z)^{n+k+1}},\
 k_{a,p}(z)=\frac{k_{a}(z)}{{\left\Vert{k_{a}}\right\Vert}_{p}}.$\ \par 
Set $\chi _{a}:={\left\Vert{k_{a}}\right\Vert}_{k,2}^{-2}=k_{a}(a)^{-2}=(1-\left\vert{a}\right\vert
 ^{2})^{n+k+1}.$ Then  $\displaystyle {\left\Vert{k_{a}}\right\Vert}_{p}\simeq
 \chi _{a}^{-1/p'}.$\ \par 
\quad And we have $\forall f\in H^{p}({\mathbb{B}}),\ f(a):={\left\langle{f,k_{a}}\right\rangle},$
 where ${\left\langle{\cdot ,\cdot }\right\rangle}$ is the scalar
 product of the Hilbert space $A_{k}^{2}({\mathbb{B}}).$\ \par 
\quad In the unit disc the $A^{p}({\mathbb{D}})$ interpolating sequences
 where characterized in the nice papers of K. Seip~\cite{Seip93},
 ~\cite{seip1}. He used densities to do it, opposite to the product
 of Gleason distances used for $H^{p}({\mathbb{D}})$ interpolating
 sequences.\ \par 

\subsection{Links between Bergman and Hardy spaces.}
\quad We shall use the "Subordination Lemma", see~\cite{subPrinAmar12}
 and earlier~\cite{amBerg78}.\ \par 
We shall write $(z,\zeta ):=(z_{1},...,z_{n},\zeta _{1},...,\zeta
 _{k+1})$ to define a point in ${\mathbb{C}}^{n+k+1}.$ Now on
 we denote ${\mathbb{B}}_{k}$ the unit ball in ${\mathbb{C}}^{k}.$\ \par 
\quad The links between $A_{k}^{p}({\mathbb{B}}_{n})$ and $H^{p}({\mathbb{B}}_{n+k+1})$
 are given by the Subordination Lemma:\ \par 
\quad \quad \quad $f\in A_{k}^{p}({\mathbb{B}}_{n})$ iff $\tilde f(z,\zeta ):=f(z)$
 is in $H^{p}({\mathbb{B}}_{n+k+1}),$ and we have ${\left\Vert{f}\right\Vert}_{k,p}\simeq
 {\left\Vert{\tilde f}\right\Vert}_{H^{p}}.$\ \par 
We also have $f(z,\zeta )\in H^{p}({\mathbb{B}}_{n+k+1})\Rightarrow
 f(z,0)\in A_{k}^{p}({\mathbb{B}}_{n})$ with ${\left\Vert{f(\cdot
 ,0)}\right\Vert}_{k,p}\leq C{\left\Vert{f}\right\Vert}_{H^{p}}.$\ \par 
\quad The reproducing kernels for $H^{p}({\mathbb{B}}_{n+k+1})$ are:\ \par 
\quad \quad \quad $\displaystyle \forall a,z\in {\mathbb{B}}_{n+k+1,\ }\tilde k_{a}(z)=\frac{1}{(1-\bar
 a\cdot z)^{n+k+1}},\ \tilde k_{a,p}(z)=\frac{\tilde k_{a}(z)}{{\left\Vert{\tilde
 k_{a}}\right\Vert}_{p}},$ with $\displaystyle {\left\Vert{\tilde
 k_{a}}\right\Vert}_{p}\simeq (1-\left\vert{a}\right\vert ^{2})^{-(n+k+1)/p'}.$\
 \par 
Hence, not surprisingly, if $a\in {\mathbb{B}}_{n},$ we get\ \par 
\quad \quad \quad $\displaystyle \forall a,z\in {\mathbb{B}}_{n},\ \tilde k_{a}(z)=k_{a}(z),\
 \tilde k_{a,p}(z)=k_{a,p}(z),$ with $\displaystyle {\left\Vert{k_{a}}\right\Vert}_{p}\simeq
 (1-\left\vert{a}\right\vert ^{2})^{-(n+k+1)/p'}.$\ \par 

\subsection{Carleson measures for Bergman spaces.}
\quad In~\cite{Hastings75} the Carleson measures for the space $A_{k}^{p}({\mathbb{B}})$
 are defined. See also~\cite{CimaMercer95},~~\cite{Zhu05},  ~\cite{AbateSaraco11}
 and ~\cite{subPrinAmar12}.\ \par 

\begin{defin}
Let $\mu $ be a Borel measure in ${\mathbb{B}}_{n}$ then it is
 called a Carleson measure for $A_{k}^{p}({\mathbb{B}}_{n})$ if:\par 
\quad \quad \quad $\displaystyle \forall p\geq 1,\ \exists C>0,\ \forall f\in A_{k}^{p}({\mathbb{B}}_{n}),\
 \int_{{\mathbb{B}}_{n}}{\left\vert{f(z)}\right\vert ^{p}\,d\mu
 (z)}\leq C{\left\Vert{f}\right\Vert}_{k,p}^{p}.$\par 
\quad The same for the usual Carleson measures for  Hardy spaces:\par 
\quad \quad \quad $\displaystyle \forall p\geq 1,\ \exists C>0,\ \forall f\in H^{p}({\mathbb{B}}_{n}),\
 \int_{{\mathbb{B}}_{n}}{\left\vert{f(z)}\right\vert ^{p}\,d\mu
 (z)}\leq C{\left\Vert{f}\right\Vert}_{H^{p}({\mathbb{B}})}^{p}.$
\end{defin}
\quad By their geometric definition~\cite{subPrinAmar12}, ~\cite{Zhu05}{\it
 , }we already know that if $\mu $ is Carleson for a $p\in \lbrack
 1,\infty \lbrack $ then it is Carleson for all $p\in \lbrack
 1,\infty \lbrack .$\ \par 
\quad Let $\mu $ be a measure in ${\mathbb{B}}_{n}$ and extend it by
 $0$ in ${\mathbb{B}}_{n+k+1}.$ Call the extended measure $\tilde
 \mu .$ We have\ \par 

\begin{prop}
~\label{e1}The measure $\mu $ is a Carleson measure for $A_{k}^{p}({\mathbb{B}}_{n})$
 iff the measure $\tilde \mu $ is a Carleson measure for $H^{p}({\mathbb{B}}_{n+k+1}).$
\end{prop}

      Proof.\ \par 
$\bullet $ Suppose first that $\tilde \mu $ is Carleson for $H^{p}({\mathbb{B}}_{n+k+1}).$\
 \par 
Then take $f\in A_{k}^{p}({\mathbb{B}}_{n}).$ Extend it as $\tilde
 f(z,\zeta ):=f(z),$ then $\tilde f\in H^{p}({\mathbb{B}}_{n+k+1}).$
 We have, because $\tilde \mu $ is supported by ${\mathbb{B}}_{n}$:\ \par 
\quad \quad \quad $\displaystyle \int_{{\mathbb{B}}_{n+k+1}}{\left\vert{\tilde
 f(z,\zeta )}\right\vert ^{p}\,d\tilde \mu (z,\zeta )}=\int_{{\mathbb{B}}_{n}}{\left\vert{f(z)}\right\vert
 ^{p}\,d\mu (z)}.$\ \par 
\quad But $\tilde \mu $ is Carleson for $H^{p}({\mathbb{B}}_{n+k+1})$
 and $\tilde f\in H^{p}({\mathbb{B}}_{n+k+1}),$ so we get\ \par 
\quad \quad \quad $\displaystyle \int_{{\mathbb{B}}_{n}}{\left\vert{f(z)}\right\vert
 ^{p}\,d\mu (z)}=\int_{{\mathbb{B}}_{n+k+1}}{\left\vert{\tilde
 f(z,\zeta )}\right\vert ^{p}\,d\tilde \mu (z,\zeta )}\leq C{\left\Vert{\tilde
 f}\right\Vert}_{H^{p}({\mathbb{B}}_{n+k+1})}^{p}\leq C{\left\Vert{f}\right\Vert}_{k,p}^{p}.$\
 \par 
Hence $\mu $ is Carleson for $A_{k}^{p}({\mathbb{B}}_{n}).$\ \par 
$\bullet $ Suppose now that $\mu $ is Carleson for $A_{k}^{p}({\mathbb{B}}_{n}).$
 Let $g(z,\zeta )\in H^{p}({\mathbb{B}}_{n+k+1}).$ We still have\ \par 
\quad \quad \quad $\displaystyle \int_{{\mathbb{B}}_{n+k+1}}{\left\vert{g(z,\zeta
 )}\right\vert ^{p}\,d\tilde \mu (z,\zeta )}=\int_{{\mathbb{B}}_{n}}{\left\vert{g(z,0)}\right\vert
 ^{p}\,d\mu (z)}.$\ \par 
Hence, because $\mu $ is Carleson for $A_{k}^{p}({\mathbb{B}}_{n}),$
 and $g(z,0)\in A_{k}^{p}({\mathbb{B}}_{n})$ by the Subordination Lemma:\ \par 
\quad \quad \quad $\displaystyle \int_{{\mathbb{B}}_{n+k+1}}{\left\vert{g(z,\zeta
 )}\right\vert ^{p}\,d\tilde \mu (z,\zeta )}=\int_{{\mathbb{B}}_{n}}{\left\vert{g(z,0)}\right\vert
 ^{p}\,d\mu (z)}\leq C{\left\Vert{g(\cdot ,0)}\right\Vert}_{k,p}^{p}\leq
 C{\left\Vert{g}\right\Vert}_{H^{p}({\mathbb{B}}_{n+k+1})}^{p}.$\ \par 
Hence $\tilde \mu $ is Carleson for $H^{p}({\mathbb{B}}_{n+k+1}).$
 $\blacksquare $\ \par 

\begin{rem}
Proposition~\ref{e1} defines extrinsically the Carleson measures
 for the Bergman spaces $A_{k}^{p}({\mathbb{B}}_{n}).$
\end{rem}
\quad As a Corollary of P. Thomas' Theorem~\cite{Thomas87} we get:\ \par 

\begin{thm}
Let $S$ be a sequence of points in ${\mathbb{B}}.$ If the sequence
 $S$ is $A_{k}^{1}({\mathbb{B}})$ interpolating then $S$ is a
 Carleson sequence for$A_{k}^{p}({\mathbb{B}}_{n}).$
\end{thm}
\quad Proof.\ \par 
The sequence $S$ is contained in the set $\zeta =0$ in ${\mathbb{B}}_{n+k+1},$
 with $z\in {\mathbb{B}}_{n}$ and $(z,\zeta )\in {\mathbb{B}}_{n+k+1}.$
 Because $S$ is $A_{k}^{1}$-interpolating we get \ \par 
\quad \quad \quad $\displaystyle \forall \lambda \in ???^{1}(S),\ \exists f\in
 A_{k}^{1}::\forall a\in S,\ {\left\langle{f,\ k_{a,\infty }}\right\rangle}=\lambda
 _{a}$ and $\displaystyle {\left\Vert{f}\right\Vert}_{k,1}\lesssim
 {\left\Vert{\lambda }\right\Vert}_{???^{1}(S)}.$\ \par 
Fix $\lambda \in ???^{1}(S)$ and $f\in A_{k}^{1}$ doing the interpolation.
 Consider $\tilde f(z,\zeta ):=f(z).$ Then, as we seen, $\tilde
 f\in H^{1}({\mathbb{B}}_{n+k+1}),\ \tilde f$ interpolates the
 sequence $\lambda \in ???^{1}(S)$ and ${\left\Vert{f}\right\Vert}_{H^{1}({\mathbb{B}}_{n+k+1})}\lesssim
 {\left\Vert{f}\right\Vert}_{k,1}\lesssim {\left\Vert{\lambda
 }\right\Vert}_{???^{1}(S)}.$ Hence we get that $S$ is $H^{1}({\mathbb{B}}_{n+k+1})$-interpolating
 in ${\mathbb{B}}_{n+k+1}.$ We apply P. Thomas' Theorem~\cite{Thomas87}
 to get that  $S$ is a Carleson sequence for $H^{p}({\mathbb{B}}_{n+k+1}).$
 Because $S\subset {\mathbb{B}}_{n}$ we deduce that $S$ is a
 Carleson sequence for$A_{k}^{p}({\mathbb{B}}_{n}).$ $\blacksquare $\ \par 

\subsection{The main result for the Bergman spaces.}

\begin{thm}
Let $S$ be a sequence of points in ${\mathbb{B}}_{n}.$ Suppose
 that $S$ is $A_{k}^{p}$-interpolating in ${\mathbb{B}}_{n}$
 for a $p\geq 2.$ then $S$ has the LEP with range $E_{S}^{p}:=\mathrm{S}\mathrm{p}\mathrm{a}\mathrm{n}\lbrace
 k_{a,p},\ a\in S\rbrace .$
\end{thm}
\quad Proof.\ \par 
The sequence $S$ is contained in the set $\zeta =0$ in ${\mathbb{B}}_{n+k+1},$
 with $z\in {\mathbb{B}}_{n}$ and $(z,\zeta )\in {\mathbb{B}}_{n+k+1}.$
 Because $S$ is $A_{k}^{p}$-interpolating we get \ \par 
\quad \quad \quad $\displaystyle \forall \lambda \in ???^{p}(S),\ \exists f\in
 A_{k}^{p}::\forall a\in S,\ {\left\langle{f,\ k_{a,p'}}\right\rangle}=\lambda
 _{a}$ and ${\left\Vert{f}\right\Vert}_{k,p}\lesssim {\left\Vert{\lambda
 }\right\Vert}_{???^{p}(S)}.$\ \par 
Fix $\lambda \in ???^{p}(S)$ and $f\in A_{k}^{p}$ doing the interpolation.
 Consider $\tilde f(z,\zeta ):=f(z).$ Then, as we seen, $\tilde
 f\in H^{p}({\mathbb{B}}_{n+k+1}),\ \tilde f$ interpolates the
 sequence $\lambda \in ???^{p}(S)$ and ${\left\Vert{f}\right\Vert}_{H^{p}({\mathbb{B}}_{n+k+1})}\lesssim
 {\left\Vert{f}\right\Vert}_{k,p}\lesssim {\left\Vert{\lambda
 }\right\Vert}_{???^{p}(S)}.$ Hence we get that $S$ is $H^{p}({\mathbb{B}}_{n+k+1})$-interpolating
 in ${\mathbb{B}}_{n+k+1}.$ So we can apply Theorem~\ref{lE17}
 which gives that $S$ has the LEP in $\mathrm{S}\mathrm{p}\mathrm{a}\mathrm{n}\lbrace
 \tilde k_{a,p},\ a\in S\rbrace .$ This gives the result because
 when $a\in {\mathbb{B}}_{n},$ we have that the reproducing kernels
 $\tilde k_{a,p}$ for $H^{p}({\mathbb{B}}_{n+k+1})$ are the same
 as the reproducing kernels $k_{a,p}$ for $A_{k}^{p}$ i.e. $\tilde
 k_{a,p}(z,\zeta )=k_{a,p}(z).$\ \par 
The proof is complete. $\blacksquare $\ \par 

\subsection{Strictly $A_{k}^{p}$ interpolating sequences.}
\quad As we did for Hardy spaces, we define the Grammian for a sequence
 of points in ${\mathbb{B}}_{n}$: $\displaystyle G_{a,b}:={\left\langle{k_{a,p'},\
 k_{b,p}}\right\rangle}.$\ \par 
Again this define an operator on $???^{p'}(S)$:\ \par 
\quad \quad \quad $\displaystyle \forall \mu \in ???^{p'}(S),\ G\mu :=\lbrace (G\mu
 )_{b}\rbrace _{b\in S},\ (G\mu )_{b}:=\sum_{a\in S}{\mu _{a}G_{a,b}}={\left\langle{\sum_{a\in
 S}{\mu _{a}k_{a,p'},k_{b,p}}}\right\rangle}.$\ \par 
We still have\ \par 

\begin{prop}
Let $S$ be a sequence of points in ${\mathbb{B}}_{n}.$ Let:\par 
\quad (CS) The sequence $S$ is Carleson for $A_{k}^{p}.$\par 
\quad (BG) The associated Grammian $G$ is bounded on $???^{p'}(S).$\par 
Then (CS) implies (BG).
\end{prop}
\quad Proof.\ \par 
We lift everything on ${\mathbb{B}}_{n+k+1}$ and, because the
 reproducing kernels agree, we get that $\tilde G_{a,b}=G_{a,b}.$
 On the other hand $S$ Carleson for $A_{k}^{p}$ means that the
 associated measure $\mu _{S}$ is Carleson for $A_{k}^{p},$ hence
 its extension $\tilde \mu _{S}$ is Carleson for $H^{p}({\mathbb{B}}_{n+k+1})$
 by Proposition~\ref{e1}. So we get, applying Proposition~\ref{lE18},
 that $\tilde G$ is bounded on $???^{p'}(S).$\ \par 
Because $\tilde G_{a,b}=G_{a,b}$ we get that $G$ is bounded on
 $???^{p'}(S).$ $\blacksquare $\ \par 

\begin{defin}
Let $S$ be a sequence of points in ${\mathbb{B}}_{n}.$ We shall
 say that $S$ is {\bf strictly }$A_{k}^{p}\ ${\bf  interpolating}
 if the operator $G$ is bounded below on $???^{p'}(S).$
\end{defin}
\quad Now we have, as for the $H^{p}$ case:\ \par 

\begin{thm}
~\label{lE21}Let $S$ be a sequence of points in ${\mathbb{B}}_{n}.$
 Suppose that $S$ is Carleson for $A_{k}^{p},$ and $S$ is strictly
 $A_{k}^{p}$ interpolating.\par 
Then $S$ is $A_{k}^{p}$ interpolating and has the LEP with range
 $E_{S}^{p}:=\mathrm{S}\mathrm{p}\mathrm{a}\mathrm{n}\lbrace
 k_{a,p},\ a\in S\rbrace .$
\end{thm}
\quad Proof.\ \par 
Lifting the situation from ${\mathbb{B}}_{n}$ to ${\mathbb{B}}_{n+k+1},$
 because the reproducing kernels agree, we get that $\tilde G_{a,b}=G_{a,b}.$
 Hence $S$ strictly $A_{k}^{p}$ interpolating implies that $S$
 is strictly $H^{p}({\mathbb{B}}_{n+k+1})$ interpolating. It
 remains to apply Theorem~\ref{lE13} to have that $S$ is $H^{p}({\mathbb{B}}_{n+k+1})$
 interpolating and has the LEP with range $\tilde E_{S}^{p}:=\mathrm{S}\mathrm{p}\mathrm{a}\mathrm{n}\lbrace
 \tilde k_{a,p},\ a\in S\rbrace .$ But $S\subset {\mathbb{B}}_{n}$
 then the reproducing kernels agree, so we get that $S$ is $A_{k}^{p}$
 interpolating and has the LEP with range $E_{S}^{p}:=\mathrm{S}\mathrm{p}\mathrm{a}\mathrm{n}\lbrace
 k_{a,p},\ a\in S\rbrace .$\ \par 
The proof is complete. $\blacksquare $\ \par 
\ \par 
\quad As for Hardy spaces we have the converse of the Theorem~\ref{lE21}.\ \par 

\begin{thm}
Let $S$ be a sequence of points in ${\mathbb{B}}_{n}.$ Suppose
 that $S$ is $A_{k}^{p}$ interpolating and has the LEP with range
 $E_{S}^{p}:=\mathrm{S}\mathrm{p}\mathrm{a}\mathrm{n}\lbrace
 k_{a,p},\ a\in S\rbrace $ for a $p>2.$\par 
Then $S$ is strictly $A_{k}^{p}$ interpolating.
\end{thm}
\quad Proof.\ \par 
For the proof we repeat exactly the arguments we use for proving
 Theorem~\ref{lE21}. $\blacksquare $\ \par 
\ \par 
\quad We also have, still using the Subordination Lemma:\ \par 

\begin{thm}
Let $S$ be a dual bounded sequence in $A_{k}^{p}.$ Then $S$ is
 $A_{k}^{s}$ interpolating with the LEP, for any $s\in \lbrack 1,p\lbrack .$
\end{thm}
\quad Proof.\ \par 
We lift the situation to ${\mathbb{B}}_{n+k+1}$ which leads to
 deal with the Hardy space $H^{p}({\mathbb{B}}_{n+k+1})$ and
 we apply Theorem 1.4 p. 482 in~\cite{intBall09}. Hence we get
 the result that $S$ is $H^{s}({\mathbb{B}}_{n+k+1})$ interpolating
 with the LEP, for any $s\in \lbrack 1,p\lbrack .$ It remains
 to go back to ${\mathbb{B}}_{n}$ to end the proof of this theorem.
 $\blacksquare $\ \par 

\section{Babenko examples in the class of holomorphic functions.}
\quad Recall the result of Babenko~\cite{Babenko47}: there exists $\lbrace
 e_{a}\rbrace _{a\in S}$ a basis of unit vectors in the Hilbert
 space ${\mathcal{H}}$ which is $2$-\emph{hilbertian} but $\lbrace
 e_{a}\rbrace _{a\in S}$ is not a Riesz basis.\ \par 
\quad The same way there exists $\lbrace e_{a}\rbrace _{a\in S}$ a
 basis of unit vectors in the Hilbert space ${\mathcal{H}}$ which
 is $2$-\emph{besselian} but $\lbrace e_{a}\rbrace _{a\in S}$
 is not a Riesz basis.\ \par 
\quad Hence a natural question is:\ \par 
\quad \quad \quad is it possible to have Babenko examples in the case of systems
 of reproducing kernels?\ \par 
\quad In the case of the unit disc the answer is no because of the
 following remarkable property of the disc algebra. Let $S$ be
 a \textbf{finite sequence} of distinct points in the unit disc.
 Let $\lbrace k_{a,2}\rbrace _{a\in S}$ be the normalized Cauchy
 kernels and set $E_{S}:=\mathrm{S}\mathrm{p}\mathrm{a}\mathrm{n}\lbrace
 k_{a,2},\ a\ \in S\rbrace .$ Set also $\lbrace \rho _{a}\rbrace
 _{a\in S}$ the dual basis of $\lbrace k_{a,2}\rbrace _{a\in
 S}$ in $E_{S}$ and $\lbrace \rho _{a,2}\rbrace _{a\in S}$ this
 normalized dual basis. \ \par 
Then we have~\cite{AmarThesis77}:\ \par 

\begin{thm}
Let $S$ be a finite sequence of points in ${\mathbb{D}}.$ The
 following anti-linear isometry is true\par 
\quad \quad \quad $\displaystyle {\left\Vert{\sum_{a\in S}{\lambda _{a}k_{a,2}}}\right\Vert}_{H^{2}({\mathbb{D}})}={\left\Vert{\sum_{a\in
 S}{\bar \lambda _{a}\rho _{a,2}}}\right\Vert}_{H^{2}({\mathbb{D}})}.$ 
\end{thm}
\quad This Theorem easily implies that if $\lbrace k_{a,2}\rbrace _{a\in
 S}$ is hilbertian or Besselian for any sequence of points in
 ${\mathbb{D}},$ then it is a Riesz sequence. Hence the Babenko
 phenomenon cannot exist in the case of the Cauchy kernels in the disc.\ \par 
\quad We have the following Theorem:\ \par 

\begin{thm}
(~\cite{Thomas87}) ~\label{lE7}Let $S$ be a sequence of points
 in ${\mathbb{B}}.$ If the sequence $S$ is $H^{1}({\mathbb{B}})$
 interpolating then $S$ is a Carleson sequence.
\end{thm}
\quad The Theorem by P. Thomas is valid also for a class of harmonic
 functions in the ball. See ~\cite{AmarWirtBoule07} for an easier
 proof using Wirtinger inequalities, but working only for holomorphic
 functions.\ \par 
\quad Using the link between interpolating sequences and besselian
 systems and between Carleson sequences and hilbertian systems,
 we have that, in $H^{p}({\mathbb{B}}),$ if a sequence of reproducing
 kernels is $p$-besselian then it is automatically a $p$-Riesz sequence.\ \par 
\quad We also have:\ \par 

\begin{thm}
(~\cite{CarlesonInter19}) Let $S$ be a sequence of points in
 ${\mathbb{D}}^{n}.$ Suppose that $S$ is $H^{p}({\mathbb{D}}^{n})$
 interpolating with a $p>2.$ Then $S$ is Carleson.
\end{thm}
\quad Then, the same way as above, we have that, in $H^{p}({\mathbb{D}}^{n}),\
 p>2,$ if a sequence of reproducing kernels is $p$-besselian
 then it is automatically a $p$-Riesz sequence.\ \par 
\quad But we have:\ \par 

\begin{thm}
The are sequences of reproducing kernels in $H^{p}(\Omega )$
 which are $p$-hilbertian but not  $p$-Riesz sequences.
\end{thm}
\quad Proof.\ \par 
By the characterisation of K. Seip~\cite{Seip93}, ~\cite{seip1},
 we know that, for any $p\geq 1,q>p,$ there are $p$-interpolating
 sequences for the Bergman space $A^{p}({\mathbb{D}})$ which
 are not $q$-interpolating for $A^{q}({\mathbb{D}}).$ \ \par 
\quad Take such a sequence $S\subset {\mathbb{D}}$ with $q>p>2.$ Then,
 using the subordination lemma, we have that $\displaystyle \tilde
 S:=\lbrace (a,0),\ a\in S\rbrace .$ can be seen as a $p$-interpolating
 sequence in $H^{p}({\mathbb{B}}_{2}),$ hence it is a $p$-Carleson
 sequence in $H^{p}({\mathbb{B}}_{2}).$ So it is also a $q$-Carleson
 sequence in $H^{q}({\mathbb{B}}_{2}).$ By Proposition~\ref{e1},
 we get that $S$ is also $q$-Carleson sequence in $A^{q}({\mathbb{D}}),$
 but it is not a $q$-interpolating sequence. Hence the normalised
 associated reproducing kernels in $\displaystyle A^{q}({\mathbb{D}})$
 make a $q$-hilbertian sequence but not a $q$-Riesz sequence.\ \par 
\quad We have the same result in $\displaystyle H^{q}({\mathbb{B}}_{2})$
 for the normalised associated reproducing kernels in $\displaystyle
 H^{q}({\mathbb{B}}_{2})$ associated to the sequence $\displaystyle
 \tilde S:=\lbrace (a,0),\ a\in S\rbrace .$\ \par 
\quad To deal with $H^{p}({\mathbb{D}}^{2})$ we set the sequence $S$
 in the diagonal of the bi-disc: $\tilde S:=\lbrace (a,a),\ a\in
 S\rbrace .$ Then again we know that $\tilde S$ is $\displaystyle
 H^{p}({\mathbb{D}}^{2})$ but not $\displaystyle H^{q}({\mathbb{D}}^{2})$
 interpolating, see~\cite{AmMen02} and the references therein.
 Because $p>2,$ we get that $\tilde S$ is $p$-Carleson, hence
 $q$-Carleson. Because $\tilde S$  is not a $q$-interpolating
 sequence, the normalised associated reproducing kernels in $\displaystyle
 H^{q}({\mathbb{D}}^{2})$ make a $q$-hilbertian sequence but
 not a $q$-Riesz sequence. $\blacksquare $\ \par 
\ \par 
This leads to the natural conjectures, still with $\Omega $ being
 the ball or the polydisc:\ \par 
\ \par 
$\bullet $ If $S$ is dual bounded in $H^{2}(\Omega )$ (resp.
 in $H^{p}(\Omega )$), then it is $H^{2}(\Omega )$ (resp. $H^{p}(\Omega
 )$) interpolating.\ \par 
\ \par 
\quad Of course these conjectures are true in one variable, because
 dual boundedness is easily seen to be equivalent to the Carleson
 condition. Hence it implies interpolation.\ \par 

\bibliographystyle{/usr/local/texlive/2017/texmf-dist/bibtex/bst/base/apalike}

\end{document}